\documentclass[a4paper,11pt]{article}
\usepackage{bbm}
\usepackage{amsfonts}
\usepackage{indentfirst,latexsym,bm,amsmath,amssymb,amsthm}
\usepackage[dvips]{graphicx}

\makeatletter\@addtoreset{equation}{section} \makeatother
 

\title{Global existence of critical nonlinear wave equation with time dependent variable coefficients}
\author{Yi Zhou
\thanks{School of Mathematical Sciences, Fudan University, Shanghai 200433, P. R. China;
Key Laboratory of Mathematics for Nonlinear Sciences (Fudan
University), Ministry of Education of China,  P. R. China; Shanghai
Key Laboratory for contemporary Applied Mathematics, School of
Mathematical Sciences, Fudan University  ({\tt
Email:yizhou@fudan.ac.cn})} \and Ning-An Lai
\thanks{School of Mathematical Sciences, Fudan University,
Shanghai 200433, P. R. China; ({\tt
Email:071018029@fudan.edu.cn})}.}
\date{}
\begin{document}
\maketitle
\begin{abstract}
In this paper, we establish global existence of smooth solutions for
the Cauchy problem of the critical nonlinear wave equation with time
dependent variable
coefficients in three space dimensions\\
\begin{equation}\nonumber\\
\partial_{tt}\phi-\partial_{x_i}\big(g^{ij}(t,x)\partial_{x_j}\phi\big)+\phi^5=0,~~\mathbb{R}_t \times
\mathbb{R}_x^3,\\
\end{equation}\\
%%where~$\Box_g$~is the wave operator associated to a given metric
%%\begin{equation}\nonumber\\
%%g=-dt^2+g_{ij}(t,x)dx^idx^j,~~~i,j=1,2,3,\\
%\end{equation}\\
where~$\big(g_{ij}(t,x)\big)$~is a regular function valued in the
spacetime of~$3\times3$~ positive definite matrix
and~$\big(g^{ij}(t,x)\big)$~its inverse matrix. Here and in the
sequence, a repeated sum on an index in lower and upper position is
never indicated. In the constant coefficients case, the result of
global existence is due to Grillakis \cite{Grillakis1}; and in the
time-independent variable coefficients case, the result of global
existence and regularity is due to Ibrahim and Majdoub
\cite{Ibrahim}. The key point of our proofs is to show that the
energy cannot concentrate at any point. For that purpose, following
Christodoulou and Klainerman \cite{Chris}, we use a null frame
associated to an optical function to construct a geometric
multiplier similar to the well-known Morawetz multiplier. Then we
use comparison theorem originated from Riemannian Geometry to
estimate the error terms. Finally, using Strichartz inequality due
to \cite{Smith} as Ibrahim and Majdoub \cite{Ibrahim}, we obtain
global existence.
\par {\bf Keywords:}  time dependent variable
 coefficients, critical nonlinearity, null frame,
 comparison theorem.
\end{abstract}
\section{\textbf{Introduction}}

\indent   In this work we consider global existence of smooth
solutions of the Cauchy problem
\begin{equation}\label{eq1}
\left \{
\begin{aligned}
&\partial_{tt}\phi-\partial_{x_i}\big(g^{ij}(t,x)\partial_{x_j}\phi\big)+\phi^5=0,~~\mathbb{R}_t
\times
\mathbb{R}_x^3,\\
&\phi(t_0,x)=f_1(x)\in
C_{0}^{\infty}(\mathbb{R}_x^3),~~~\phi_{t}(t_0,x)=f_2(x)\in
C_{0}^{\infty}(\mathbb{R}_x^3),\\
\end{aligned} \right.
\end{equation}
here~$\{g^{ij}(t,x)\}_{i,j=1}^3$~denotes a matrix valued smooth
function
 of the variables~$(t,x)\in \mathbb{R} \times \mathbb{R}^3$,~
 which takes values in the real, symmetric,~$3\times3$~matrices,
 such that for some~$C >0$,~
\begin{equation}\label{eq101}
C|\xi|^2\leq g^{ij}(t, x)\xi_i\xi_j \leq
C^{-1}|\xi|^2,~~~~\forall~\xi \in
\mathbb{R}^3,~(t, x)\in \mathbb{R} \times \mathbb{R}^3 .\\
\end{equation}
Obviously it is a critical wave equation on a curved spacetime.
First let us survey existence and regularity results for critical
nonlinear wave equations briefly. If ~$g^{ij}=\delta^{ij}$,~which
denotes the Kronecker delta function, we say
problem~$\eqref{eq1}$~is of constant coefficients. In the case of
critical nonlinear wave equation with constant coefficients, a
wealth of results are available in the literature. For cauchy
problem,  global existence of~$C^2$-solutions in dimension~$n=3$~was
first obtained by Rauch \cite{Rauch}, assuming the initial energy to
be small. In 1988, also for "large" data global~$C^2$-solutions in
dimension~$n=3$~were shown to exist by Struwe \cite{Struwe} in the
radially symmetric case. Grillakis \cite{Grillakis1} in 1990 was
able to remove the latter symmetry assumption and obtained the same
result. Not much later, Kapitanskii \cite{Ka} estiblished the
existence of a unique, partially regular solution for all
dimensions. Combining Strichartz inequality and Morawetz estimates,
Grillakis \cite{Grillakis2} in 1992 established global existence and
regularity for dimensions~$3\leq n\leq 5$~and announced the
corresponding results in the radial caes for dimensions~$n\leq 7$.~
Then Shatah and Struwe \cite{Shatah} obtained global existence and
regularity for dimensions~$3\leq n\leq 7$.~They also proved the
global well-posedness in the energy space in \cite{Shatah2} 1994.
For the critical exterior problem in dimension 3, Smith and Sogge
\cite{Sogge} in 1995 proved global existence of smooth solutions. In
2008, Burq et all \cite{Burq} obtained the same result in 3-D
bounded domain. \\
\indent For the critical Cauchy problem with time-independent
variable coefficients, Ibrahim and Majdoub \cite{Ibrahim} in 2003
studied the existence of both global smooth for dimensions~$3\leq n<
6$~and Shatah-Struwe's solutions for dimensions~$n\geq 3$.~Recently,
we have showed global existence and regularity in \cite{Zhou} for
the critical exterior problem with time-independent
variable coefficients in dimension~$n=3$.~\\
\indent In this paper we are interested in the critical case with
coefficients depending on the time and space variables. Our result
concerns global existence and regularity, showed as follow: \\
{\bf Theorem 1.1.}  Problem~$\eqref{eq1}$~admits a unique global
solution~$\phi \in C^\infty(\mathbb{R}\times \mathbb{R}^3)$.~\\
\indent The demonstration of theorem 1.1 is done by contradiction,
showing~$\phi$~is uniformly bounded. For that purpose, the key step
is to show the non-concentration of the~$L^6$~part of the energy
(and hence the energy), and to do this the idea is to work in
geodesic cone just like
light cone in constant coefficients case. Thus we have\\
{\bf Lemma 1.2.}~(\textbf{Non-concentration lemma}) If~$\phi \in
([t_0, 0)\times \mathbb{R}^3)$~solves \eqref{eq1}, then
\begin{equation}\nonumber\\
\lim_{t\rightarrow 0}\int_{Q(t)}\phi^6\mathrm{d}v=0,\\
\end{equation}
where~$Q(t)$~is the intersection of~$t$~time slice with backward
solid
characteristic cone from origin.\\
\indent In order to prove the non-concentration lemma, in the
constant coefficients caes the Morawetz
multiplier~$t\partial_t+r\partial_r+1$~is used, where~$r=|x|$;~while
in the time-independent variable coefficients case the geometric
multiplier~$t\partial_t+\rho\nabla \rho+1$~is used instead,
where~$\rho$~is the associated distance function. The time-dependent
variable coefficients case considered in this work is much more
difficult, and the simple minded generalization to use
multiplier~$t\partial_t+(\underline{u}-t)\nabla(\underline{u}-t)+1$~will
not work, where~$\underline{u}$~is an optical function (close
to~$t+|x|$~). Following Christodoulou and Klainerman \cite{Chris} we
use a null frame. However, the emphasis in their work is the
asymptotic behavior of the null frame at infinity, and here we
emphasize its asymptotic behavior locally at a possible blow up
point. We derive the asymptotic behavior of the null frame by using
comparison theorem originated
from Riemainnian geometry. \\
\indent To prove our result, we also need Strichartz inequality, stated as \\
{\bf Lemma 1.3.}~(\textbf{Strichartz inequality})
 Assuming~$g^{ij}(t,x)$~satisfies the conditions of the
introduction,~$\phi$~solves the Cauchy problem as follow in the half
open strip~$[t_0, 0)\times \mathbb{R}^3$:~
\begin{equation}\nonumber\\
\left \{
\begin{aligned}
&\partial_{tt}\phi-\partial_{x_i}\big(g^{ij}(t,x)\partial_{x_j}\phi\big)=F(t,x),\\
&\phi(t_0,x)=f_1(x)\in
C_{0}^{\infty}(\mathbb{R}_x^3),~~~\phi_{t}(t_0,x)=f_2(x)\in
C_{0}^{\infty}(\mathbb{R}_x^3),\\
\end{aligned} \right.
\end{equation}
then we have
\begin{eqnarray}\label{104}
\|\phi\|_{L_{t}^{\frac{2q}{q-6}}L_{x}^q([t_0,~0)\times~\mathbb{R}^3)}\leq
C\big(\|f_1\|_{H^1(\mathbb{R}^3)}+\|f_2\|_{L^2(\mathbb{R}^3)}+
\|F\|_{L_{t}^1L_{x}^2([t_0,~0)\times~\mathbb{R}^3)}\big)\nonumber \\
6\leq q <\infty.
 \end{eqnarray}
 \indent For the proof see Smith \cite{Smith}.\\
 \indent Then combining these two lemmas we
 can
 establish the uniform bounds on the local solution~$\phi$,~which
 implies our result, this step is completely
 parallel to Ibrahim and Majdoub \cite{Ibrahim} and we omit it.\\
 \indent Our results can extend to a more general variable
 coefficient second order partial differential equation with operator as follow:\\
\begin{equation}\nonumber\\
\begin{aligned}
L\equiv
\partial_t^2+2b^i(t,x)\partial_{ti}^2-a^{ij}(t,x)\partial_{ij}^2+L_1,~~a^{ij}=a^{ji},\\
\end{aligned}
\end{equation}
where all coefficients are real and~$C^\infty$,~and~$L_1$~is a first
order operator. In general, we can get rid of cross terms (that is,
terms like~$
b^i\partial_{ti}^2$~) by the following procedure: let us write (with new first order terms~$L_1'$~)\\
\begin{equation}\nonumber\\
\begin{aligned}
L\equiv
(\partial_t+b^i\partial_i)^2-(a^{ij}+b^ib^j)\partial_{ij}^2+L_1'.\\
\end{aligned}
\end{equation}
If, in the region under consideration, we can perform a change of
variables
\begin{equation}\nonumber\\
\begin{aligned}
X_1=\varphi_1(t,x),\cdots, X_n=\varphi_n(t,x), T=t,\\
\end{aligned}
\end{equation}
and set
\begin{equation}\nonumber\\
\begin{aligned}
\frac{\partial \varphi_j}{\partial t}+b^i\frac{\partial
\varphi_j}{\partial x_i}=0,~~j=1,\cdots,n,\\
\end{aligned}
\end{equation}
 in such a way that the vector
field~$\partial_t+b^i\partial_i$~becomes~$\partial_T$,~then the
operator~$L$~takes the form
\begin{equation}\nonumber\\
\begin{aligned}
\overline{L}=
\partial_T^2-\overline{a}^{ij}\partial_{X_iX_j}^2+\overline{L}_1,\\
\end{aligned}
\end{equation}
for some new coefficients~$\overline{a}^{ij}$~and lower order
terms~$\overline{L}_1$.~\\
\indent As an application of our result, we can consider the
critical wave equation in the Schwarzschild spacetime~$(\mathcal
{M}, g)$~with parameter~$M > 0$,~where~$g$~is the Schwarzschild
metric whose line
element is\\
\begin{equation}\nonumber\\
\begin{aligned}
ds^2=-\Big(1-\frac{2M}{r}\Big)dt^2+\Big(1-\frac{2M}{r}\Big)^{-1}dr^2+r^2d\omega^2,\\
\end{aligned}
\end{equation}
where~$d\omega^2$~is the measure on the
sphere~$\mathbbm{s}^2$.~While the singularity at~$r=0$~is a true
metric singularity, we note that the apparent singularity
at~$r=2M$~is merely a coordinate singularity. Indeed, define the
Regge-Wheeler tortoise
coordinate~$r^*$~by\\
\begin{equation}\nonumber\\
\begin{aligned}
r^*=r+2M\log(r-2M)-3M-2M\log M,\\
\end{aligned}
\end{equation}
and set~$v=t+r^*$.~Then in the~$(r^*, t, \omega)$~coordinates the
Schwarzschild metric~$g$~
is expressed in the form\\
\begin{equation}\nonumber\\
\begin{aligned}
ds^2=-\Big(1-\frac{2M}{r}\Big)dv^2+2dvdr+r^2d\omega^2.\\
\end{aligned}
\end{equation}
Let~$\Sigma$~be an arbitrary Cauchy surface for the (maximally
extended) Schwarzschild spacetime~$(\mathcal {M},
g)$~stated as above and consider the Cauchy problem of the wave equation\\
\begin{equation}\label{1000}\\
\left \{
\begin{aligned}
&\Box_g \phi-\phi^5=0,\\
&(\phi ,\phi_t )|_{\Sigma}=(\psi_0,\psi_1), \\
\end{aligned} \right.
\end{equation}
for this problem Marzuola et all \cite{Mar} proved global existence
and uniqueness of finite energy solution under the assumption of
small initial energy, and
according to our result we can remove the small energy assumption, that is\\
{\bf Theorem 1.4.} For smooth initial data prescribed on~$\Sigma$,
equation \eqref{1000} admits a unique
global smooth solution in the~$(r^*, t, \omega)$~coordinates.\\
\indent Now we sketch the plan of this article. In the next
section
 we recall some geometric concepts which are necessary for our proofs.
 Section 3 is devoted to the proof of lemma 1.2: the fundamental lemma
 expressing the non-concentration of
 ~$L^6$~part of the energy.\\
 \indent Finally, we remark that although our non-concentration
 lemma is stated only in dimension~$n=3$,~the proof works in any
 dimension for the critical wave equations.\\
 \indent In this paper, the letter~$C$~denotes a constant which may
change from one to the other.\\
\section{\textbf{Null frame}}
\indent Let~$\{g_{ij}(t,x)\}_{i,j=1}^3$~denotes the inverse matrix
of~$\{g^{ij}(t,x)\}_{i,j=1}^3$,~and consider the spilt
metric~$g=-dt^2+g_{ij}(t,x)dx^idx^j=g_{\alpha\beta}dx^\alpha
dx^\beta$~on~$\mathbb{R}_{x,t}^4$~(close to the Minkowski metric).
So we will work in the spacetime, a 4-dimensional
manifold~$M.$~Local coordinates on~$M$~are denoted by~$x_\alpha,
\alpha=1,2,3,4.$~The convention is used that Latin indices run from
1 to 3 while Greek indices relate to the spacetime manifold~$M$~and
run from 1 to 4. The index 4 corresponds to the time coordinate,
while~$(x_1,x_2,x_3)$~are the spatial coordinates. The corresponding
partial derivatives are~$\partial_\alpha=\frac{\partial}{\partial
x_\alpha}.$~We introduce an optical function~$\underline{u}$~(close
to~$t+|x|$~)for~$g$:~ a~$C^1$~function which satisfies the eikonal
equation\\
\begin{equation}\label{eq105}
\left \{
\begin{aligned}
&g^{\alpha\beta}\partial_\alpha\underline{u}\partial_\beta\underline{u}
=g_{\alpha\beta}\partial^\alpha\underline{u}\partial^\beta\underline{u}
=\langle\nabla\underline{u},
\nabla\underline{u}\rangle=|\nabla\underline{u}|^2=0,\\
&\underline{u}(t,0)=t,\\
\end{aligned} \right.
\end{equation}
where~$\langle ~,~ \rangle$~denotes the inner product about the
given metric. In PDE terms, this means that the level
surfaces~$\{\underline{u}=C\}$~are characteristic surfaces for any
operator with principal
symbol~$g^{\alpha\beta}\xi_\alpha\xi_\beta$.~From this construction
it is easy to see that the first
order derivatives of~$\underline{u}$~are locally bounded.\\
\indent Then we set
 \begin{equation}\label{eq0}
 \begin{aligned}
&\underline{L}=-\nabla
\underline{u}=(\partial_t\underline{u})\partial_t-(g^{ij}\partial_i\underline{u})\partial_j=
m^{-1}(\partial_t+N),\\
&L=\frac{\partial_t}{\partial_t\underline{u}}+\frac
{(g^{ij}\partial_i\underline{u})\partial_j}{(\partial_t\underline{u})^2}=m(\partial_t-N),\\
\end{aligned}
\end{equation}\\
where~$\nabla$~is the gradient about the given
metric,~$m=(\partial_t\underline{u})^{-1},$~~$N=-\frac{(g^{ij}\partial_i\underline{u})\partial_j}
{\partial_t\underline{u}}=-(mg^{ij}\partial_i\underline{u})\partial_j.$
~It is easily to see that they are close
to~$\partial_t-\partial_r$~and~$\partial_t+\partial_r$~respectively.
And~$D_{\underline{L}}\underline{L}=0$,~showing that a integral
curve of~$\underline{L}$~is a geodesic. This follows from the
symmetry of the Hessian, since for any vector field~$X$,~we have
 \begin{equation}\nonumber\\
 \begin{aligned}
&<D_{\underline{L}}\underline{L},X>=-<D_{\underline{L}}\nabla\underline{u},X>
=-<D_X\nabla\underline{u},\underline{L}>\\
&=<D_X\underline{L},\underline{L}>=\frac{1}{2}X<\underline{L},\underline{L}>
=\frac{1}{2}X<\nabla\underline{u},\nabla\underline{u}>=0.\\
\end{aligned}
\end{equation}\\
So the integral curves of the field~$\underline{L}$~generate a
backward geodesic cone with vertices on the~$t$-axis.
~%~Then using Hessian and
%Laplace comparison theorems from Lorentaian Geometry and Riemannian
%Geometry to estimate the error terms we show the~$L^6$~part of the
%energy associated to (1.1) cannot concentrate at any point~$(t_0,
%x_0)$,~which is the key step to get global existence.\\
%
%
%
%In this section we show the main result and proofs.\\
%Following Alinhac [], let~$g=-dt^2+g_{ij}(t,x)dx_idx_j$~be a split
%metric on~$\mathbb{R}_{x,t}^4$~(close to the Minknowski metric)
%and~$\underline{u}$~an optical function for~$g$~(close to~$t+r$~).
%So~$\underline{u}$~is a~$C^1$~function which satisfies the eikonal
%equation\\
%\begin{equation}\nonumber\\
%\left \{
%\begin{aligned}
%&g^{\alpha\beta}\partial_\alpha\underline{u}\partial_\beta\underline{u}
%=g_{\alpha\beta}\partial^\alpha\underline{u}\partial^\beta\underline{u}
%=\langle\nabla\underline{u},
%\nabla\underline{u}\rangle=|\nabla\underline{u}|^2=0,\\
%&\underline{u}(t,0)=t,\\
%\end{aligned} \right.
%\end{equation}
%\begin{equation}\nonumber\\
%g^{\alpha\beta}\partial_\alpha\underline{u}\partial_\beta\underline{u}
%=g_{\alpha\beta}\partial^\alpha\underline{u}\partial^\beta\underline{u}
%=\langle\triangledown\underline{u},
%\triangledown\underline{u}\rangle=|\triangledown\underline{u}|^2=0.\\
%\end{equation}
 Using the coordinate~$t$,~we define the
foliation~$\sum_{t_1}=\{(x,t),t=t_1\}$,~and
using~$\underline{u}$,~we
define the foliation by nonstandard 2-spheres as\\
\begin{equation}\nonumber\\
S_{t_1,\underline{u}_1}=\{(x,t),t=t_1,
\underline{u}(x,t)=\underline{u}_1\}.\\
\end{equation}
Since~$\triangledown\underline{u}$~ is orthogonal
to~$\{\underline{u}=\underline{u}_1\}$~and~$\partial_t$~is
orthogonal to~$\sum_{t_1}$,~the field~$\underline{L}$~is a null
vector orthogonal to the geodesic cone and~$N$~is an horizontal
field
orthogonal to~$S_{t_1,\underline{u}_1}$.~Moreover,\\
\begin{equation}\nonumber\\
\langle N,
N\rangle=\frac{1}{(\partial_t\underline{u})^2}g_{ij}(g^{ik}\partial_k\underline{u})(g^{jl}\partial_l
\underline{u})=\frac{1}{(\partial_t\underline{u})^2}g^{kl}\partial_k\underline{u}\partial_l\underline{u}
=1.\\
\end{equation}
Then, if~$(e_1, e_2)$~form an orthonormal basis on the nonstandard
spheres, the frame\\
\begin{equation}\nonumber\\
e_1, e_2, e_3\equiv L=m(\partial_t-N), e_4\equiv
\underline{L}=m^{-1}(\partial_t+N)\\
\end{equation}
is a null frame with
\begin{equation}\nonumber\\
\left \{
\begin{aligned}
&\langle e_1, e_1\rangle=\langle e_2, e_2\rangle=1, \langle e_1,
e_2\rangle=0,\\
&\langle e_1, L\rangle=\langle e_1, \underline{L}\rangle=\langle
e_2, L\rangle=\langle e_2, \underline{L}\rangle=0,\\
&\langle L, L\rangle=\langle \underline{L}, \underline{L}\rangle=0,
\langle L, \underline{L}\rangle=-2.\\
\end{aligned} \right.
\end{equation}
\indent We will work in the null frame as above and it requires that
we know the vector~$D_{\alpha}e_{\beta}$,~that is: the frame
coefficients~$<D_{\alpha}e_{\beta},e_{\gamma}>$.~\\
\indent We define the frame coefficients by
\begin{equation}\label{eq600}
\begin{aligned}
&<D_a\underline{L},e_b>=\underline{\chi}_{ab}=\underline{\chi}_{ba},~~
<D_aL,e_b>=\chi_{ab}=\chi_{ba},\\
&<D_{\underline{L}}\underline{L},e_a>=0, ~~~~~~~~~~~~~<D_LL,e_a>=2\xi_a,\\
&<D_{\underline{L}}L,e_a>=2\eta_a,
~~~~~~~~~~<D_L\underline{L},e_a>=2\underline{\eta}_a,\\
&<D_{\underline{L}}\underline{L},L>=0,
~~~~~~~~~~~~~~<D_LL,\underline{L}>=4\underline{\omega}=-<D_L\underline{L},L>,\\
\end{aligned}
\end{equation}
where~$a,b=1,2$.~\\
%Setting~$X=\frac{1}{2}(L+\underline{L})$,~one can compute the
%components of the deformation tensor~$^{(X)}\!\pi=\pi$:~\\
%\begin{equation}\nonumber\\
%\begin{aligned}
%&\pi_{\underline{L}\underline{L}}=0, \pi_{\underline{L}L}=2\omega,
%\pi_{LL}=-4\underline{\omega},\\
%&\pi_{\underline{L}e_a}=\eta_a-\underline{\eta}_a,
%\pi_{Le_a}=\xi_a+2\underline{\eta}_a,
%\pi_{ab}=\underline{\chi}_{ab}+\chi_{ab}.\\
%\end{aligned}
%\end{equation}
\indent If we call~$k$~the second fundamental form of~$\Sigma_t$~by
\begin{equation}\nonumber\\
k(X,Y)=-<D_X\partial_t,Y>, k_{ij}=-\frac{1}{2}\partial_tg_{ij},\\
 \end{equation}
then~$k$~is nothing but the first order derivatives of~$g$~and so
bounded. By some simple computation, we also have
\begin{equation}\label{eq116}
\begin{aligned}
2\eta_a&=-2k_{Na},\\
2\underline{\eta}_a&=-2me_a(\underline{u}_t)+2k_{Na},\\
2\xi_a&=-2m^2\underline{\eta}_a+2m^2k_{Na},\\
\chi_{ab}&=-m^2\underline{\chi}_{ab}-2mk_{ab},\\
\underline{\omega}&=-\partial_tm=\frac{\partial_{tt}\underline{u}}{(\partial_t\underline{u})^2}.\\
\end{aligned}
 \end{equation}
For the details, one can read Alinhac \cite{Alinhac}. And we are
interested in the asymptotic behavior of these frame coefficients
near the origin, thus we have\\
% {\bf Lemma 2.1.} Assuming~$\xi_a,
%~\eta_a, ~\underline{\eta}_a$~are
%frame coefficients as above, then they are all bounded locally. \\
%{\bf Proof.}  As~$[e_a,\partial_t](\underline{u})$~is nothing but
%the first order derivative of~$\underline{u}$,~all the first order
%derivatives of the optical function~$\underline{u}$~and~$k$~are
%bounded locally as mentioned above, then the lemma follows from
%\eqref{eq116}.\\
{\bf Theorem 2.1.} Assuming~$\xi_a, ~\eta_a, ~\underline{\eta}_a,~
\underline{\omega},~ \underline{\chi}_{ab}, ~\chi_{ab}$~are frame
coefficients as above,
then\\
\begin{align}
&|\eta_a|\leq C,\\
 &\frac{ct}{2}\leq
\underline{\omega}=\frac{\partial_{tt}\underline{u}}{(\partial_t\underline{u})^2}
=m^2\partial_{tt}\underline{u}\leq -\frac{ct}{2},\\
&\frac{1}{t-\underline{u}}+ct\leq \underline{\chi}_{aa}\leq
\frac{1}{t-\underline{u}}-ct,\\
&|\underline{\eta}_{a}|\leq -Ct,\\
&|\xi_a|\leq C-Ct,\\
&|1-m||\underline{\chi}_{aa}|\leq C,\\
&4+Ct\leq \chi_{aa}\underline{u}+\chi_{bb}\underline{u}+
\underline{\chi}_{aa}u+\underline{\chi}_{bb}u\leq 4-Ct,\\
&(2+Ct)|\overline{\nabla}\phi|^2\leq
\sum_{a,b=1}^{2}(\chi_{ab}\underline{u}+\underline{\chi}_{ab}u)e_a(\phi)e_b(\phi)
\leq (2-Ct)|\overline{\nabla}\phi|^2,
\end{align}
where~$a=1, 2;~c,~C$~are positive constants;~$t< 0$~as we work in
the backward geodesic cone starting from the origin
and~$u=2t-\underline{u}$.~\\
\indent The first inequality is obviously from \eqref{eq116}, and to
prove the other inequality of this theorem, we need several
lemmas stated below. First let us introduce the comparison theorem.\\
\indent Assuming~$C,D$~take values in the real,
symmetric,~$(n-1)\times(n-1)$~matrices. If for
any~$(\alpha_1,\cdots,\alpha_{n-1}),(\beta_1,\cdots,\beta_{n-1})\in
\mathbb{R}^{n-1}$~and~$\sum\alpha_i^2=\sum\beta_i^2$,~we have
\begin{equation}\nonumber\\
\left(
                                \begin{array}{ccc}
                                  \alpha_1 , \cdots , \alpha_{n-1} \\
                                \end{array}
                              \right)
C\left(
                                 \begin{array}{c}
                                   \alpha_1 \\
                                   \vdots \\
                                   \alpha_{n-1} \\
                                 \end{array}
                               \right)\geq \left(
                                \begin{array}{ccc}
                                  \beta_1 , \cdots , \beta_{n-1} \\
                                \end{array}
                              \right)
D\left(
                                 \begin{array}{c}
                                   \beta_1 \\
                                   \vdots \\
                                   \beta_{n-1} \\
                                 \end{array}
                               \right),\\
\\
\end{equation}
then we say~$C\succ D$.~\\
 {\bf Lemma 2.2.} Let~$gl(n-1, \mathbb{R})$~be set of~$n-1$~order
 real symmetric matrix,~$K,~\widetilde{K}:[0, b)\rightarrow gl(n-1,
 \mathbb{R})$.~Suppose~$A:[0, b)\rightarrow gl(n-1,
 \mathbb{R})$~satisfies the ordinary differential equation
\begin{equation}\nonumber\\
\left \{
\begin{aligned}
&\frac{d^2A}{ds^2}+AK=0,\\
&A(0)=0, ~\frac{dA}{ds}(0)=I (the ~unit ~matrix),\\
\end{aligned} \right.
\end{equation}
and~$\widetilde{A}:[0, b)\rightarrow gl(n-1,
 \mathbb{R})$~satisfies
\begin{equation}\nonumber\\
\left \{
\begin{aligned}
&\frac{d^2\widetilde{A}}{ds^2}+A\widetilde{K}=0,\\
&\widetilde{A}(0)=0, ~\frac{d\widetilde{A}}{ds}(0)=I (the ~unit ~matrix),\\
\end{aligned} \right.
\end{equation}
where~$s\in [0, b)$.~Also we assume~$A,\widetilde{A}$~are invertible
in~$[0, b)$~and~$K\prec \widetilde{K}$,~then
\begin{equation}\label{eq106}
A^{-1}\frac{dA}{ds}\succ \widetilde{A}^{-1}\frac{\widetilde{A}}{ds}.\\
 \end{equation}
\indent For the proof see \cite{Wu}.\\
\indent If the
metric~$\widetilde{g}=-dt^2+\widetilde{g}_{ij}(x)dx^idx^j$,~$\widetilde{g}_{ij(x)}$~depend
only on the spatial coordinates,
then~$\widetilde{\underline{u}}=t+\widetilde{\rho}$~is an optical
function
for~$\widetilde{g}$~satisfies~$\widetilde{\underline{u}}(t,0)=t$,~where~$\widetilde{\rho}$~is
the Riemannian distance function on the Riemannian
manifold~$(\mathbb{R}^3, ~\widetilde{g}_{ij}(x))$.~The corresponding
null frame related to~$\widetilde{\underline{u}}$~is
\begin{equation}\nonumber\\
\begin{aligned}
\widetilde{e}_1,~ \widetilde{e}_2,
~\widetilde{e}_3=\partial_t+\partial_{\widetilde{\rho}}=\partial_t+\widetilde{g}^{ij}\widetilde{\rho}_i\partial_j,
~\widetilde{e}_4=\partial_t-\partial_{\widetilde{\rho}}=\partial_t-\widetilde{g}^{ij}\widetilde{\rho}_i\partial_j,\\
\end{aligned}
\end{equation}
where~$\{\widetilde{g}^{ij}(x)\}_{i,j=1}^3$~denotes the inverse
matrix of~$\{\widetilde{g}_{ij}(x)\}_{i,j=1}^3$.~And
then~\underline{u}~can be compared
with~$\widetilde{\underline{u}}$~through lemma 2.2.\\
\indent Let~$\gamma: [0, b)\rightarrow (\mathbb{R}_{x,t}^4, g)$~is
the integral curve of~$\underline{L}=-\nabla u$~and we call it null
geodesic, as~$<\underline{L},
\underline{L}>=0$,~then~$\dot{\gamma}=\underline{L}=-\nabla
u$.~Let~$\{e_1, e_2, e_3=L, e_4=\underline{L}\}$~be parallel null
frame along~$\gamma$,~$J_i(s)$~be the Jacobi field
along~$\gamma$,~satisfies~$J_i(0)=0, \dot{J_i}(0)=e_i(0),
(i=1,2,3)$.~So we have
\begin{equation}\nonumber\\
\left(
  \begin{array}{c}
    J_1(s) \\
    J_2(s) \\
    J_3(s) \\
  \end{array}
  \right)=A(s)\left(
                    \begin{array}{c}
                      e_1(s) \\
                      e_2(s) \\
                      e_3(s) \\
                    \end{array}
                  \right),\\
 \end{equation}
where~$A(s)$~denotes an invertible matrix valued function of the
parameter~$s\in [0, b)$.~Then the Jacobi equation becomes
\begin{equation}\nonumber\\
\frac{d^2A}{ds^2}+AK=0,\\
 \end{equation}
where~$K=(K_{ij})_{i,j=1}^3, K_{ij}=<R(\dot{\gamma},
e_i)\dot{\gamma}, e_j>$~denotes a symmetric matrix of~$3\times3$.~We
then easily get\\
\begin{equation}\label{eq109}
H_{\underline{u}}(e_i, e_j)=D^2\underline{u}(e_i, e_j)=-(A^{-1}\frac{dA}{ds})_{ij}=-\underline{\chi}_{ij},\\
\end{equation}
where~$H_{\underline{u}}$~denotes the Hessian form
of~$\underline{u}$.~And then \eqref{eq600} yields
\begin{equation}\label{eq601}
(\underline{\chi}_{ij})_{i,j=1}^3=\left(
                                    \begin{array}{ccc}
                                      \underline{\chi}_{11} & \underline{\chi}_{12} & 2\underline{\eta}_1 \\
                                      \underline{\chi}_{12} & \underline{\chi}_{22} & 2\underline{\eta}_2 \\
                                      2\underline{\eta}_1 & 2\underline{\eta}_2 & -4\underline{\omega} \\
                                    \end{array}
                                  \right).\\
 \end{equation}
Correspondingly for optical
function~$\widetilde{\underline{u}}=t+\widetilde{\rho}$,~we
have~$\widetilde{A}, \widetilde{K}$~and
\begin{equation}\label{eq151}
H_{\widetilde{\underline{u}}}(\widetilde{e}_i,
\widetilde{e}_j)=D^2\widetilde{\underline{u}}(\widetilde{e}_i,
\widetilde{e}_j)=D^2(t+\widetilde{\rho})(\widetilde{e}_i,
\widetilde{e}_j)=
-(\widetilde{A}^{-1}\frac{d\widetilde{A}}{ds})_{ij}=-\widetilde{\underline{\chi}}_{ij}.\\
 \end{equation}
\indent Note that one assumption of lemma 2.2 is~$K\prec
\widetilde{K}$,~but~$K_{33}=<R(\underline{L},~ L)\underline{L},~
L>\neq0$,~while
~$\widetilde{K}_{33}=<R(\partial_t-\partial_{\widetilde{\rho}},~
\partial_t+\partial_{\widetilde{\rho}})\partial_t-\partial_{\widetilde{\rho}},~
\partial_t+\partial_{\widetilde{\rho}}>=0$,~so we have to introduce a
conformally related metric tensor to~$\widetilde{g}$~to ensure the
condition~$K\prec \widetilde{K}$.~Let~$(\widetilde{g}_{ij}(x),
\mathbb{R}^3)$~be a space form with positive constant sectional
curvature~$c$.~We set then the conformally related
metric~$\widetilde{g}_c=e^{ct^2}\widetilde{g}$.~\\
{\bf Lemma 2.3.} Let~(M, g)~be a semi-Riemannian manifold of
dimension~$n$~and let~$g_c=\varphi g$~be a conformally related
metric tensor to~$g$,~where~$\varphi: M\rightarrow (0, \infty)$~is a
map. Then \\
(1) ~$\mathop \nabla
\limits^c=\frac{1}{\varphi}\nabla$,~where~$\nabla$~and~$\mathop
\nabla \limits^c$~are the gradients on~$(M, g)$~and~$(M,
g_c)$,~respectively.\\
(2) ~For~$X, Y\in \Gamma TM$,~
\begin{equation}\nonumber\\
{\mathop
\nabla\limits^c}_XY=\nabla_XY+\frac{1}{2\varphi}X(\varphi)Y+\frac{1}{2\varphi}Y(\varphi)X
-\frac{1}{2\varphi}g(X, Y)\nabla \varphi,\\
 \end{equation}
where~$\nabla$~and~$\mathop \nabla \limits^c$~are the Levi-Civita
connections of~$(M,g)$~and~$(M,
g_c)$,~respectively.\\
(3) ~If~$f: M\rightarrow \mathbb{R}$~is a map then, for~$X, Y\in
\Gamma TM$,\\
\begin{equation}\nonumber\\
 \begin{aligned}
H_f^c(X, Y)=&H_f(X, Y)-\frac{1}{2\varphi}[g(\nabla\varphi,
X)g(\nabla f, Y)\\
&+g(\nabla f, X)g(\nabla\varphi, Y)-g(\nabla\varphi, \nabla
f)g(X, Y)],\\
 \end{aligned}
 \end{equation}
 where~$H_f$~and~$H_f^c$~are the Hessian forms of~$f$~on~$(M,g)$~and~$(M,
g_c)$,~respectively.\\
\indent For the proof see \cite{Eduardo}.\\
\indent From lemma 2.3, we have
\begin{equation}\nonumber\\
 \begin{aligned}
<\mathop \nabla\limits^c(t+\widetilde{\rho}), \mathop
\nabla\limits^c(t+\widetilde{\rho})>_{\widetilde{g}_c}=e^{ct^2}<\frac{1}{e^{ct^2}}\nabla(t+\widetilde{\rho}),
\frac{1}{e^{ct^2}}\nabla(t+\widetilde{\rho})>_{\widetilde{g}}=0,\\
 \end{aligned}
 \end{equation}
then~$t+\widetilde{\rho}$~is also an optical function
for~$\widetilde{g}_c$.~As above, we define~$\widetilde{A}_c,
\widetilde{K}_c,
\widetilde{\chi}_{cij}$~associated to~$\widetilde{g}_c$.~\\
 \indent It is easily known that for a manifold~$(M, g_M)$~ with constant
 curvature~$c$~
\begin{equation}\label{eq405}
 \begin{aligned}
&R(X, Y, Z, W)\\
&=c\big[g_M(X, Z)g_M(Y, W)-g_M(X, W)g_M(Y,
Z)\big],~\forall~ X, Y, W, Z\in \Gamma TM.\\
 \end{aligned}
 \end{equation}
Then for the space form~$(\widetilde{g}_{ij}(x), \mathbb{R}^3)$~with
positive constant sectional curvature~$c$~a computation according to
\eqref{eq405} gives
\begin{equation}\nonumber\\
 \begin{aligned}
\widetilde{K}&=\left(
                \begin{array}{ccc}
                  c & 0 & 0 \\
                  0 & c & 0 \\
                  0 & 0 & 0 \\
                \end{array}
              \right),\\
\widetilde{\underline{\chi}}_{11}&=-H_{\widetilde{\underline{u}}}(\widetilde{e}_1,~
\widetilde{e}_1)=-<D_{\widetilde{e}_1}\nabla(t+\widetilde{\rho}),~
\widetilde{e}_1>_{\widetilde{g}}=-<D_{\widetilde{e}_1}\nabla t,~
\widetilde{e}_1>_{\widetilde{g}}-<D_{\widetilde{e}_1}\nabla_g\widetilde{\rho},~
\widetilde{e}_1>_{\widetilde{g}}\\
&=-D^2\widetilde{\rho}(\widetilde{e}_1,~
\widetilde{e}_1)=-\sqrt{c}\cot(\sqrt{c}\widetilde{\rho})<\widetilde{e}_1,~
\widetilde{e}_1>_{\widetilde{g}}=-\sqrt{c}\cot(\sqrt{c}\widetilde{\rho}),\\
\widetilde{\underline{\chi}}_{22}&=-H_{\widetilde{\underline{u}}}(\widetilde{e}_2,~
\widetilde{e}_2)=-\sqrt{c}\cot(\sqrt{c}\widetilde{\rho})<\widetilde{e}_2,~
\widetilde{e}_2>_{\widetilde{g}}=-\sqrt{c}\cot(\sqrt{c}\widetilde{\rho}),\\
\widetilde{\underline{\chi}}_{12}&=-H_{\widetilde{\underline{u}}}(\widetilde{e}_1,~
\widetilde{e}_2)=-\sqrt{c}\cot(\sqrt{c}\widetilde{\rho})<\widetilde{e}_1,~
\widetilde{e}_2>_{\widetilde{g}}=0,\\
\widetilde{\underline{\chi}}_{33}&=-H_{\widetilde{\underline{u}}}(\widetilde{e}_3,~
\widetilde{e}_3)=-<D_{\widetilde{e}_3}\nabla(t+\widetilde{\rho}),~
\widetilde{e}_3>_{\widetilde{g}}=<D_{\partial_t+\partial\widetilde{\rho}}\partial_t-\partial\widetilde{\rho},~
\partial_t+\partial\widetilde{\rho}>_{\widetilde{g}}=0,\\
\widetilde{\underline{\chi}}_{13}&=-H_{\widetilde{\underline{u}}}(\widetilde{e}_1,~
\widetilde{e}_3)=-<D_{\widetilde{e}_3}\nabla(t+\widetilde{\rho}),~
\widetilde{e}_1>_{\widetilde{g}}=<D_{\partial_t+\partial\widetilde{\rho}}\partial_t-\partial\widetilde{\rho},~
\widetilde{e}_1>_{\widetilde{g}}=0,\\
\widetilde{\underline{\chi}}_{23}&=-H_{\widetilde{\underline{u}}}(\widetilde{e}_2,~
\widetilde{e}_3)=-<D_{\widetilde{e}_3}\nabla(t+\widetilde{\rho}),~
\widetilde{e}_2>_{\widetilde{g}}=<D_{\partial_t+\partial\widetilde{\rho}}\partial_t-\partial\widetilde{\rho},~
\widetilde{e}_2>_{\widetilde{g}}=0,\\
 \end{aligned}
 \end{equation}
where~$\nabla_g$~is the gradient on the space
form~$(\widetilde{g}_{ij}(x), \mathbb{R}^3)$.~Thanks to lemma 2.3,
after the conformal change of metric they become
\begin{equation}\label{eq113}
 \begin{aligned}
\widetilde{K}_c&=\left(
                \begin{array}{ccc}
                  2ce^{ct^2}-c^2t^2e^{ct^2} & 0 & 0 \\
                  0 & 2ce^{ct^2}-c^2t^2e^{ct^2} & 0 \\
                  0 & 0 & 4ce^{ct^2} \\
                \end{array}
              \right),\\
\widetilde{\underline{\chi}}_{c11}&=-H_{\widetilde{\underline{u}}}^c(\widetilde{e}_1,~
\widetilde{e}_1)=-\big(H_{\widetilde{\underline{u}}}(\widetilde{e}_1,~
\widetilde{e}_1)-ct\big)=-\sqrt{c}\cot(\sqrt{c}\widetilde{\rho})+ct,\\
\widetilde{\underline{\chi}}_{c22}&=-H_{\widetilde{\underline{u}}}^c(\widetilde{e}_2,~
\widetilde{e}_2)=-\big(H_{\widetilde{\underline{u}}}(\widetilde{e}_2,~
\widetilde{e}_2)-ct\big)=-\sqrt{c}\cot(\sqrt{c}\widetilde{\rho})+ct,\\
\widetilde{\underline{\chi}}_{c12}&=-H_{\widetilde{\underline{u}}}^c(\widetilde{e}_1,~
\widetilde{e}_2)=-H_{\widetilde{\underline{u}}}(\widetilde{e}_1,~
\widetilde{e}_2)=0,\\
\widetilde{\underline{\chi}}_{c33}&=-H_{\widetilde{\underline{u}}}^c(\widetilde{e}_3,~
\widetilde{e}_3)=-\big(H_{\widetilde{\underline{u}}}(\widetilde{e}_3,~
\widetilde{e}_3)-2ct\big)=2ct,\\
\widetilde{\underline{\chi}}_{c13}&=-H_{\widetilde{\underline{u}}}^c(\widetilde{e}_1,~
\widetilde{e}_3)=-H_{\widetilde{\underline{u}}}(\widetilde{e}_1,~
\widetilde{e}_3)=0,\\
\widetilde{\underline{\chi}}_{c23}&=-H_{\widetilde{\underline{u}}}^c(\widetilde{e}_2,~
\widetilde{e}_3)=-H_{\widetilde{\underline{u}}}(\widetilde{e}_2,~
\widetilde{e}_3)=0.\\
 \end{aligned}
 \end{equation}
So we can ensure~$K\prec \widetilde{K}_c$~when~$t\rightarrow
0$~and~$c$~big enough.\\
\indent Similarly, if we let~$(\widetilde{\widetilde{g}}_{ij}(x),
\mathbb{R}^3)$~be a space form with negative constant sectional
curvature~$-c$~and set the conformally related
metric~$\widetilde{\widetilde{g}}_c=e^{-ct^2}\widetilde{g}$,~then we
have
\begin{equation}\label{eq128}
 \begin{aligned}
\widetilde{\widetilde{K}}_c&=\left(
                \begin{array}{ccc}
                  -2ce^{-ct^2}-c^2t^2e^{-ct^2} & 0 & 0 \\
                  0 & -2ce^{-ct^2}-c^2t^2e^{-ct^2} & 0 \\
                  0 & 0 & -4ce^{-ct^2} \\
                \end{array}
              \right),\\
\widetilde{\widetilde{\underline{\chi}}}_{c11}&=-H_{\widetilde{\widetilde{\underline{u}}}}^c(\widetilde{\widetilde{e}}_1,~
\widetilde{\widetilde{e}}_1)=-\big(H_{\widetilde{\widetilde{\underline{u}}}}(\widetilde{\widetilde{e}}_1,~
\widetilde{\widetilde{e}}_1)+ct\big)=-\sqrt{c}\coth(\sqrt{c}\widetilde{\widetilde{\rho}})-ct,\\
\widetilde{\widetilde{\underline{\chi}}}_{c22}&=-H_{\widetilde{\widetilde{\underline{u}}}}^c(\widetilde{\widetilde{e}}_2,~
\widetilde{\widetilde{e}_2})=-\big(H_{\widetilde{\widetilde{\underline{u}}}}(\widetilde{\widetilde{e}}_2,~
\widetilde{\widetilde{e}}_2)+ct\big)=-\sqrt{c}\coth(\sqrt{c}\widetilde{\widetilde{\rho}})-ct,\\
\widetilde{\widetilde{\underline{\chi}}}_{c12}&=-H_{\widetilde{\widetilde{\underline{u}}}}^c(\widetilde{\widetilde{e}}_1,~
\widetilde{\widetilde{e}}_2)=-H_{\widetilde{\widetilde{\underline{u}}}}(\widetilde{\widetilde{e}}_1,~
\widetilde{\widetilde{e}}_2)=0,\\
\widetilde{\widetilde{\underline{\chi}}}_{c33}&=-H_{\widetilde{\widetilde{\underline{u}}}}^c(\widetilde{\widetilde{e}}_3,~
\widetilde{\widetilde{e}}_3)=-\big(H_{\widetilde{\widetilde{\underline{u}}}}(\widetilde{\widetilde{e}}_3,~
\widetilde{\widetilde{e}}_3)+2ct\big)=-2ct,\\
\widetilde{\widetilde{\underline{\chi}}}_{c13}&=-H_{\widetilde{\widetilde{\underline{u}}}}^c(\widetilde{\widetilde{e}}_1,~
\widetilde{\widetilde{e}}_3)=-H_{\widetilde{\widetilde{\underline{u}}}}(\widetilde{\widetilde{e}}_1,~
\widetilde{\widetilde{e}}_3)=0,\\
\widetilde{\widetilde{\underline{\chi}}}_{c23}&=-H_{\widetilde{\widetilde{\underline{u}}}}^c(\widetilde{\widetilde{e}}_2,~
\widetilde{\widetilde{e}}_3)=-H_{\widetilde{\widetilde{\underline{u}}}}(\widetilde{\widetilde{e}}_2,~
\widetilde{\widetilde{e}}_3)=0.\\
 \end{aligned}
 \end{equation}
Again we can ensure~$\widetilde{\widetilde{K}}_c\prec K$.~\\
\indent Using lemma 2.2, we get
\begin{equation}\label{eq108}
A^{-1}\frac{dA}{ds}\succ \widetilde{A}_c^{-1}\frac{d\widetilde{A}_c}{ds},\\
 \end{equation}
together with~$\eqref{eq109}$,~we have
\begin{equation}\label{eq111}
\big(H_{\underline{u}}(e_i, e_j\big )\prec
\big(H_{\widetilde{\underline{u}}}^c(\widetilde{e}_i,
\widetilde{e}_j)\big).\\
 \end{equation}
So
\begin{equation}\label{eq112}
4\underline{\omega}=-<D_L\underline{L}, L>=H_{\underline{u}}(e_3,
e_3)\leq H_{\widetilde{\underline{u}}}^c(\widetilde{e}_3,
\widetilde{e}_3),\\
 \end{equation}
combining \eqref{eq113} and \eqref{eq112}, we obtain
\begin{equation}\label{eq114}
\underline{\omega}\leq -\frac{ct}{2}.\\
 \end{equation}
For the same reason we have
\begin{equation}\label{eq115}
\underline{\omega}\geq \frac{ct}{2},~\\
 \end{equation}
 and the inequality (2.6) of theorem 2.1 follows.\\
\indent Combining lemma 2.2, \eqref{eq113} and \eqref{eq128}, we
obtain
\begin{equation}\label{eq186}
 \begin{aligned}
-\sqrt{c}\cot(\sqrt{c}\widetilde{\rho})+ct=\widetilde{\underline{\chi}}_{caa}
\leq \underline{\chi}_{aa}\leq
\widetilde{\widetilde{\underline{\chi}}}_{caa}=
-\sqrt{c}\coth(\sqrt{c}\widetilde{\widetilde{\rho}})-ct,~~a=1,2,\\
 \end{aligned}
 \end{equation}
as we adopt comparison theorem along integral curves
of~$\underline{L}=-\nabla
\underline{u}$~,~$\widetilde{\underline{L}}=-\nabla
\widetilde{\underline{u}}$~and~$\widetilde{\widetilde{\underline{L}}}=-\nabla
\widetilde{\widetilde{\underline{u}}}$~we
set~$\underline{u}=\widetilde{\underline{u}}=t+\widetilde{\rho}=\widetilde{\widetilde{\underline{u}}}
=t+\widetilde{\widetilde{\rho}}$,~so when~$t,
~\widetilde{\rho},~\widetilde{\widetilde{\rho}}$~are small (close to
0), we have
\begin{equation}\label{eq612}
 \begin{aligned}
\frac{1}{t-\underline{u}}+ct\leq \underline{\chi}_{aa}\leq
\frac{1}{t-\underline{u}}-ct
,~~a=1,2,\\
 \end{aligned}
 \end{equation}
%As we adopt comparison theorem along integral curves
%of~$\underline{L}=-\nabla
%\underline{u}$~and~$\widetilde{\underline{L}}=-\nabla
%\widetilde{\underline{u}}$,~we
%set~$\underline{u}=\widetilde{\underline{u}}=t+\rho$,~that is
%\begin{equation}\label{eq130}
% \begin{aligned}
%\underline{u}-t=\rho,\\
% \end{aligned}
% \end{equation}
%and then
%\begin{equation}\label{eq185}
% \begin{aligned}
%&\underline{\chi}_{aa}=-\frac{1}{\rho}+Ct=\frac{1}{t-\underline{u}}+Ct,~~a=1,2.\\
% \end{aligned}
% \end{equation}
which is the desired inequality (2.7).\\
\indent Using lemma 2.2, together with \eqref{eq109}, we have
\begin{equation}\nonumber\\
\big(\widetilde{\underline{\chi}}_{cij}\big)_{i,j=1}^3 \prec
\big(\underline{\chi}_{ij}\big)_{i,j=1}^3\prec
\big(\widetilde{\widetilde{\underline{\chi}}}_{cij}\big)_{i,j=1}^3,\\
\end{equation}
so
\begin{equation}\nonumber\\
\begin{aligned}
&\left(
  \begin{array}{ccc}
    1 & 0 & 1 \\
  \end{array}
\right)\big(\widetilde{\underline{\chi}}_{cij}\big)\left(
                                                                \begin{array}{c}
                                                                  1 \\
                                                                  0 \\
                                                                  1 \\
                                                                \end{array}
                                                              \right)\\
                                                              &\leq \left(
  \begin{array}{ccc}
    1 & 0 & 1 \\
  \end{array}
\right)\big(\underline{\chi}_{ij}\big)\left(
                                                                \begin{array}{c}
                                                                  1 \\
                                                                  0 \\
                                                                  1 \\
                                                                \end{array}
                                                              \right)\\
                                                              &\leq \left(
  \begin{array}{ccc}
    1 & 0 & 1 \\
  \end{array}
\right)\big(\widetilde{\widetilde{\underline{\chi}}}_{cij}\big)\left(
                                                                \begin{array}{c}
                                                                  1 \\
                                                                  0 \\
                                                                  1 \\
                                                                \end{array}
                                                              \right),\\
\end{aligned}
\end{equation}
thus
\begin{equation}\label{eq609}
 \begin{aligned}
\widetilde{\underline{\chi}}_{c11}+2\widetilde{\underline{\chi}}_{c13}+\widetilde{\underline{\chi}}_{c33}
\leq
\underline{\chi}_{11}+2\underline{\chi}_{13}+\underline{\chi}_{33}\leq
\widetilde{\widetilde{\chi}}_{c11}+2\widetilde{\widetilde{\chi}}_{c13}+\widetilde{\widetilde{\chi}}_{c33}.\\
 \end{aligned}
 \end{equation}
 As from lemma 2.2 we have
\begin{equation}\nonumber\\
 \begin{aligned}
\widetilde{\underline{\chi}}_{cii}\leq \underline{\chi}_{ii}\leq
\widetilde{\widetilde{\chi}}_{cii},~~i=1,2,3,
 \end{aligned}
 \end{equation}
then \eqref{eq609} means
\begin{equation}\label{eq610}
 \begin{aligned}
\widetilde{\underline{\chi}}_{c11}+2\widetilde{\underline{\chi}}_{c13}+\widetilde{\underline{\chi}}_{c33}
-\widetilde{\widetilde{\chi}}_{c11}-\widetilde{\widetilde{\chi}}_{c33}\leq
2\underline{\chi}_{13}\leq
\widetilde{\widetilde{\chi}}_{c11}+2\widetilde{\widetilde{\chi}}_{c13}+\widetilde{\widetilde{\chi}}_{c33}
-\widetilde{\underline{\chi}}_{c11}-\widetilde{\underline{\chi}}_{c33}.\\
 \end{aligned}
 \end{equation}
Combining \eqref{eq601}, \eqref{eq113}, \eqref{eq128} and
\eqref{eq612}, a calculation gives
\begin{equation}\nonumber\\
 \begin{aligned}
\big|\underline{\eta}_1\big|\leq -Ct.
 \end{aligned}
 \end{equation}
For~$\big|\underline{\eta}_2\big|$~we have the same result, and then we obtain the inequality (2.8) in theorem 2.1.\\
\indent After that, inequality (2.9) can be easily obtained from
\eqref{eq116}.\\
\indent To prove the inequality (2.10) of theorem 2.1 we need the following two steps:\\
\indent First, show~$\mathop \nabla\limits^0\underline{u}_t$~is
bounded, where~$\mathop \nabla\limits^0$~is the gradient on
Euclidean space. From equation \eqref{eq105}, we have
\begin{equation}\nonumber\\
 \begin{aligned}
&\partial_t(g^{ij}\partial_i\underline{u}\partial_j\underline{u})=
\partial_tg^{ij}\partial_i\underline{u}\partial_j\underline{u}+2g^{ij}
\partial_i\underline{u}\partial_t\partial_j\underline{u}\\
&=\partial_tg^{ij}\partial_i\underline{u}\partial_j\underline{u}+2g^{ij}
\partial_i\underline{u}\partial_j\underline{u}_t\\
&=\partial_t(\partial_t\underline{u})^2\\
&=2\partial_t\underline{u}\partial_{tt}\underline{u}.\\
 \end{aligned}
 \end{equation}
As the first order derivatives of~$\underline{u}$~are bounded,
together with (2.6), we have
\begin{equation}\nonumber\\
 \begin{aligned}
|N(\underline{u}_t)|\leq C.\\
 \end{aligned}
 \end{equation}
 Also from \eqref{eq116}, we get
 \begin{equation}\nonumber\\
 \begin{aligned}
me_a(\underline{u}_t)=-\underline{\eta}_a+k_{Na},\\
 \end{aligned}
 \end{equation}
 thanks to inequality (2.8), it implies
\begin{equation}\nonumber\\
 \begin{aligned}
|e_a(\underline{u}_t)|\leq C,~~a=1,2,\\
 \end{aligned}
 \end{equation}
 so we finish the first step.\\
\indent Second, as the result of the first step\\
\begin{equation}\label{eq125}
 \begin{aligned}
|1-m|=\big|\frac{\partial_t{\underline{u}}-1}{\partial_t{\underline{u}}}\big|
=|m||\underline{u}_t(t, x)-\underline{u}_t(t, 0)|\leq
C\sup_x|\mathop
\nabla\limits^0\underline{u}_t||x|\leq C|x|.\\
 \end{aligned}
 \end{equation}
Set~$\underline{v}=t+\delta|x|,~ \delta>0$,~then if we
choose~$\delta$~small enough
\begin{equation}\nonumber\\
 \begin{aligned}
\dot{\gamma}(\underline{v})=\underline{L}(\underline{v})=\partial_t{\underline{u}}
\partial_t{\underline{v}}-g^{ij}\partial_i{\underline{u}}\partial_i{\underline{v}}
=\partial_t{\underline{u}}-\delta
g^{ij}\partial_i{\underline{u}}\frac{x_j}{|x|}\geq
1-C\delta >0,\\
 \end{aligned}
 \end{equation}
while
\begin{equation}\nonumber\\
 \begin{aligned}
\dot{\gamma}(\underline{u})=\underline{L}(\underline{u})=0.\\
 \end{aligned}
 \end{equation}
 As~$\gamma$~is a backwards integral curve of~$\underline{L}$,~along the curve~$\gamma$~we
 conclude
\begin{equation}\nonumber\\
 \begin{aligned}
\underline{u}\geq \underline{v}=t+\delta|x|,\\
 \end{aligned}
 \end{equation}
 thus
 \begin{equation}\label{eq126}
 \begin{aligned}
\underline{u}-t\geq \delta|x|.\\
 \end{aligned}
 \end{equation}
Combining \eqref{eq125}, \eqref{eq126} and (2.7), we have
\begin{equation}\nonumber\\
 \begin{aligned}
&|1-m||\underline{\chi}_{aa}|\leq C_1,\\
 \end{aligned}
 \end{equation}
which means inequality (2.10).\\
 {\bf Lemma 2.4.} Inside the geodesic
cone where~$\underline{u}\leq 0$~, we have
\begin{equation}\label{eq196}
\begin{aligned}
|\underline{u}|\leq C|t|,~~|u|\leq C|t|.\\
\end{aligned}
\end{equation}
{\bf Proof.} By \eqref{eq126}, along the integral curve
of~$\underline{L}$~starting from the origin, we have
\begin{equation}\label{eq305}
\begin{aligned}
t< t+\delta|x|\leq \underline{u}\leq 0,\\
\end{aligned}
\end{equation}
then
\begin{equation}\nonumber\\
\begin{aligned}
2t\leq u=2t-\underline{u}\leq t,\\
\end{aligned}
\end{equation}
so we complete the proof of lemma 2.4.\\
\indent By \eqref{eq116}, we have
\begin{equation}\nonumber\\
\begin{aligned}
&\chi_{aa}\underline{u}+\chi_{bb}\underline{u}+
\underline{\chi}_{aa}u+\underline{\chi}_{bb}u \\
&=2(\underline{\chi}_{aa}+\underline{\chi}_{bb})(t-\underline{u})+
(1-m^2)(\underline{\chi}_{aa}+\underline{\chi}_{bb})\underline{u}-2m(k_{aa}+k_{bb})\underline{u},\\
\end{aligned}
\end{equation}
then (2.7), (2.10) and lemma 2.4 yield the inequality (2.11).\\
\indent Now we prove the last inequality of theorem 2.1. Using lemma
2.2 again, we have
\begin{equation}\nonumber\\
\begin{aligned}
&\left(
  \begin{array}{ccc}
    e_1(\phi) & e_2(\phi) & 0 \\
  \end{array}
\right)\big(\widetilde{\underline{\chi}}_{cij}\big)\left(
                                                                \begin{array}{c}
                                                                  e_1(\phi) \\
                                                                  e_2(\phi) \\
                                                                  0 \\
                                                                \end{array}
                                                              \right)\\
                                                              &\leq \left(
  \begin{array}{ccc}
    e_1(\phi) & e_2(\phi) & 0 \\
  \end{array}
\right)\big(\underline{\chi}_{ij}\big)\left(
                                                                \begin{array}{c}
                                                                  e_1(\phi) \\
                                                                  e_2(\phi) \\
                                                                  0 \\
                                                                \end{array}
                                                              \right)\\
                                                              &\leq \left(
  \begin{array}{ccc}
    e_1(\phi) & e_2(\phi) & 0 \\
  \end{array}
\right)\big(\widetilde{\widetilde{\underline{\chi}}}_{cij}\big)\left(
                                                                \begin{array}{c}
                                                                  e_1(\phi) \\
                                                                  e_2(\phi) \\
                                                                  0 \\
                                                                \end{array}
                                                              \right),\\
\end{aligned}
\end{equation}
together with \eqref{eq113} and \eqref{eq128}, we arrive at
\begin{equation}\nonumber\\
\begin{aligned}
&\big(-\sqrt{c}\cot(\sqrt{c}\widetilde{\rho})+ct\big)\big((e_1(\phi))^2+(e_1(\phi))^2\big)\\
&\leq
\sum_{a,b=1}^{2}\underline{\chi}_{ab}e_a(\phi)e_b(\phi)\\
&\leq
\big(-\sqrt{c}\coth(\sqrt{c}\widetilde{\widetilde{\rho}})-ct\big)\big((e_1(\phi))^2+(e_1(\phi))^2\big),\\
\end{aligned}
\end{equation}
which implies (~$t, \widetilde{\rho}, \widetilde{\widetilde{\rho}}
$~small)
\begin{equation}\label{eq350}
\begin{aligned}
\big(\frac{1}{t-\underline{u}}+ct\big)|\overline{\nabla}\phi|^2\leq
\sum_{a,b=1}^{2}\underline{\chi}_{ab}e_a(\phi)e_b(\phi)\leq
\big(\frac{1}{t-\underline{u}}-ct\big)|\overline{\nabla}\phi|^2,\\
\end{aligned}
\end{equation}
and \eqref{eq116} yields
\begin{equation}\label{eq155}
\begin{aligned}
&\sum_{a,b=1}^{2}(\chi_{ab}\underline{u}+\underline{\chi}_{ab}u)e_a(\phi)e_b(\phi)\\
&=2(t-\underline{u})\sum_{a,b=1}^{2}\underline{\chi}_{ab}e_a(\phi)e_b(\phi)+(1-m^2)\underline{u}
\sum_{a,b=1}^{2}\underline{\chi}_{ab}e_a(\phi)e_b(\phi)\\
&-2m\underline{u}\sum_{a,b=1}^{2}k_{ab}e_a(\phi)e_b(\phi).\\
\end{aligned}
\end{equation}
As~$k_{ab}$~is bounded, combining \eqref{eq350}, (2.10) and lemma
2.4 we conclude
\begin{equation}\label{eq351}
\begin{aligned}
(2+Ct)|\overline{\nabla}\phi|^2\leq \sum_{a,b=1}^{2}(\chi_{ab}\underline{u}+\underline{\chi}_{ab}u)e_a(\phi)e_b(\phi)
\leq(2-Ct)|\overline{\nabla}\phi|^2.\\
\end{aligned}
\end{equation}
So we finish the proof of theorem 2.1.\\
\section{\textbf{Non-concentration of the~$L^6$~part of the energy}}
\indent In this section, we will prove lemma 1.2, which is essential
to prove global existence and regularity. First we introduce some notations.\\
\indent Let~$z_0=(0, 0)$,~be the vertices of the backward geodesic
cone, then
\begin{equation}\nonumber\\
 \begin{aligned}
Q(z_0)=\{(t, x)\in [t_0, 0)\times \mathbb{R}^3:~\underline{u}\leq
0,~~t_0<0\},\\
 \end{aligned}
 \end{equation}
denotes the backward geodesic cone, if~$t_0\leq s_1 < s_2 <0$,~set
\begin{equation}\nonumber\\
 \begin{aligned}
Q_{s_1}^{s_2}=Q(z_0)\cap([s_1, s_2]),\\
 \end{aligned}
 \end{equation}
and
\begin{equation}\nonumber\\
 \begin{aligned}
M_{s_1}^{s_2}=\partial Q_{s_1}^{s_2}=\{(t, x)\in Q_{s_1}^{s_2}:~\underline{u}=0\},\\
 \end{aligned}
 \end{equation}
denotes the mantle associated with the truncated cone~$Q_{s_1}^{s_2}$.~\\
\begin{equation}\nonumber\\
 \begin{aligned}
Q(s)=\{x\in \mathbb{R}^3:~\underline{u}\leq 0,~t=s\}\\
 \end{aligned}
 \end{equation}
denotes the spatial cross-sections of the backward
cone~$Q(z_0)$~when
the time is~$s$.~\\
\indent Define the energy of problem~$\eqref{eq1}$~
\begin{equation}\label{eq2}
E_1(t)=\frac{1}{2}\int_{\mathbb{R}^3}
\Big(\phi_{t}^2+g^{ij}(t,x)\partial_i\phi\partial_j\phi+\frac{\phi^6}{3}\Big)\,\mathrm{d}x.\\
\end{equation}
As we have showed in section 2 that~$\partial_{tt}\underline{u},
\mathop \nabla\limits^0\underline{u}_t$~are bounded locally,
then~$\underline{u}_t$~is continuous and together with \eqref{eq105}
we have
\begin{equation}\label{eq131}
\lim_{t, x\rightarrow 0}m(t,x)=\frac{1}{\lim_{t, x\rightarrow
0}\partial_t\underline{u}(t,x)}=\frac{1}{\lim_{t, x\rightarrow
0}\partial_t\underline{u}(t,0)}=1,\\
\end{equation}
that is:~$m=1+\mathcal {O}(t).$~So when~$t$~is small,~$E_1(t)$~has a
equivalent form
\begin{equation}\label{eq3}
E(t)=\frac{1}{4}\int_{\mathbb{R}^3}
\Big(m^{-1}(L(\phi))^2+m(\underline{L}(\phi))^2
+(m+m^{-1})|\overline{\nabla}\phi|^2+\frac{m+m^{-1}}{3}\phi^6\Big)\,\mathrm{d}v,\\
 \end{equation}
where~$|\overline{\nabla}\phi|^2=\big(e_1(\phi)\big)^2+\big(e_2(\phi)\big)^2,$~and~$dv=\sqrt{|g|}dx$~is
the volume element corresponding to the metric~$g$.~Denoting the
energy density
\begin{equation}\nonumber\\
e(t)=\frac{1}{4}\Big(m^{-1}(L(\phi))^2+m(\underline{L}(\phi))^2
+(m+m^{-1})|\overline{\nabla}\phi|^2+\frac{m+m^{-1}}{3}\phi^6\Big).\\
 \end{equation}
 \indent We then define the energy flux across~$M_s^t$:~
 \begin{equation}\label{eq320}
Flux_1(\phi,
M_s^t)=\int_{M_s^t}\frac{\frac{\partial_t\underline{u}}{2}\big(\phi_t^2+g^{ij}\partial_i\phi
\partial_j\phi+\frac{\phi^6}{3}\big)-\phi_tg^{ij}\partial_i\underline{u}\partial_j\phi}
{\sqrt{(\partial_t\underline{u})^2+\sum_{j=1}^3(
\partial_j\underline{u})^2}}\mathrm{d}\nu,\\
\end{equation}
where~$d\nu$~denotes the induced Lebesgue measure
on~$M_s^t.$~Similar to the energy, it has an equivalent form
when~$t$~is small
\begin{equation}\label{eq152}
Flux(\phi,
M_s^t)=\int_{M_s^t}\frac{|\overline{\nabla}\phi|^2+\big(\underline{L}(\phi)\big)^2+\frac{\phi^6}{3}}
{2\sqrt{(\partial_t\underline{u})^2+(g^{ij}\partial_i\underline{u})^2}}\mathrm{d}\sigma,\\
\end{equation}
where~$d\sigma=\sqrt{|g|}d\nu$~denotes the volume element
corresponding to the metric~$g$~on~$M_S^T$,~and it implies
\begin{equation}\nonumber\\
Flux(\phi,
M_s^t)\geq 0.\\
\end{equation}
 {\bf Lemma 3.1.} When~$t$~is small,~$E_1(t)$~and~$Flux_1(\phi,
M_s^t)$~are equivalent
 to~$E(t)$~and~$Flux(\phi,
M_s^t)$~respectively, that is:~$E_1(t)\backsimeq E(t), Flux_1(\phi,
M_s^t)\backsimeq Flux(\phi,
M_s^t)$.~\\
Proof. Since
\begin{equation}\nonumber\\
\underline{L}= m^{-1}(\partial_t+N), L=m(\partial_t-N),\\
\end{equation}
we get
\begin{equation}\label{eq500}
\partial_t=\frac{1}{2}(m^{-1}L+m\underline{L}),\\
\end{equation}
so
 \begin{equation}
 \begin{aligned}
(\partial_t\phi)^2 &=\big[\frac{1}{2}\big(m^{-1}L(\phi)
+m\underline{L}(\phi)\big)\big]^2\\
&=\frac{1}{4}\big(m^{-2}\big(L(\phi)\big)^2\big)+2L(\phi)\underline{L}(\phi)
+m^2\big(\underline{L}(\phi)\big)^2\big).\nonumber\\
\end{aligned}
\end{equation}\\
And
\begin{equation}\nonumber\\
\langle\nabla\phi,\nabla\phi\rangle=-(\partial_t\phi)^2+g^{ij}\partial_i\phi\partial_j\phi
=\big(e_1(\phi)\big)^2+\big(e_2(\phi)\big)^2-L(\phi)\underline{L}(\phi),\\
\end{equation}
which yield
\begin{equation}\label{eq175}
g^{ij}\partial_i\phi\partial_j\phi=|\overline{\nabla}\phi|^2-L(\phi)\underline{L}(\phi)+
(\partial_t\phi)^2,\\
\end{equation}
then we get
\begin{equation}\nonumber\\
\begin{aligned}
&E_1(t)=\frac{1}{4}\int_{\mathbb{R}^3}
\Big(m^{-2}(L(\phi))^2+m^2(\underline{L}(\phi))^2
+2|\overline{\nabla}\phi|^2+\frac{2}{3}\phi^6\Big)\,\mathrm{d}x,\\
&Flux_1(\phi,
M_s^t)=\int_{M_s^t}\frac{1}{\sqrt{(\partial_t\underline{u})^2+\sum_{j=1}^3(
\partial_j\underline{u})^2}}\Big[\frac{1}{2\partial_t\underline{u}}\big((\partial_t\underline{u}\phi_t)^2
-2\partial_t\underline{u}\phi_tg^{ij}\partial_i\underline{u}\partial_j\phi\\
&+ (g^{ij}\partial_i\underline{u}\partial_j\phi)^2\big)
+\frac{\partial_t\underline{u}}{2}\big(g^{ij}\partial_i\phi\partial_j\phi+\frac{\phi^6}
{3}\big)-\frac{1}{2\partial_t\underline{u}}\big(g^{ij}\partial_i\underline{u}\partial_j\phi\big)^2\Big]\mathrm{d}\nu\\
&=\int_{M_s^t}\frac{1}{\sqrt{(\partial_t\underline{u})^2+\sum_{j=1}^3(
\partial_j\underline{u})^2}}\Big[\frac{m}{2}\big(\underline{L}(\phi)\big)^2+\frac{1}{2m}\big(
|\overline{\nabla}\phi|^2-L(\phi)\underline{L}(\phi)+
(\partial_t\phi)^2\big)\\
&+\frac{\phi^6}{6m}
-\frac{\partial_t\underline{u}}{2}\Big(\frac{g^{ij}\partial_i\underline{u}\partial_j\phi}{\partial_t\underline{u}
}\Big)^2\Big]\mathrm{d}\nu\\
&=\int_{M_s^t}\frac{1}{\sqrt{(\partial_t\underline{u})^2+\sum_{j=1}^3(
\partial_j\underline{u})^2}}\Big[\frac{m}{2}\big(\underline{L}(\phi)\big)^2+\frac{1}{2m}\big(
|\overline{\nabla}\phi|^2\\
&-m^{-1}(\partial_t\phi+N(\phi))m(\partial_t\phi-N(\phi))+
(\partial_t\phi)^2\big) +\frac{\phi^6}{6m}
-\frac{1}{2m}\big(N(\phi)\big)^2\Big]\mathrm{d}\nu\\
&=\int_{M_s^t}\frac{
\frac{1}{m}|\overline{\nabla}\phi|^2+m\big(\underline{L}(\phi)\big)^2
+\frac{1}{3m}\phi^6}{2\sqrt{(\partial_t\underline{u})^2+\sum_{j=1}^3(
\partial_j\underline{u})^2}}\mathrm{d}\nu,\\
\end{aligned}
\end{equation}
together with~$\eqref{eq131}$,~we obtain the result.\\
\indent To finish the proof, We shall require several other lemmas.
The first is standard and says that the energy associated with our
equation is
bounded.\\
{\bf Lemma 3.2.} If ~$\phi\in C^\infty([t_0, 0)\times
\mathbb{R}^3)$~is a solution to (1.1), then~$E_1(t)$~or~$E(t)$~is
bounded for all~$t_0\leq t <0$.~Additionally, if ~$t_0\leq s< t
<0$,~then
 \begin{equation}\label{eq135}
Flux(\phi, M_s^t)\rightarrow 0,~~when~s, t\rightarrow 0.\\
 \end{equation}
{\bf Proof.} To prove the boundedness of energy one multiplies both
sides of the equation
$\phi_{tt}-\frac{\mathrm{\partial}}{\mathrm{\partial}x_{i}}(g^{ij}(t,x)\phi_j)+\phi^5=0$
by $\partial_{t}\phi$ to obtain the identity
 \begin{equation}\label{eq138}
\frac{\partial}{\partial{t}}\Big(\frac{\phi_{t}^2+g^{ij}(t,x)\phi_i\phi_j}{2}+\frac{\phi^6}{6}\Big)
-\frac{1}{2}\partial_tg^{ij}(t,x)\phi_i\phi_j-\frac{\mathrm{\partial}}
{\mathrm{\partial}x_{i}}\Big(\phi_{t}g^{ij}(t,x)\phi_j\Big)=0.
\end{equation}\\
Thus,
 \begin{equation}\label{eq140}
\begin{aligned}
&\frac{\partial}{\partial{t}}\int_{\mathbb{R}^3}\Big(\frac{\phi_{t}^2+g^{ij}(t,x)
\phi_i\phi_j}{2}+\frac{\phi^6}{6}\Big)\mathrm{d}x
-\int_{\mathbb{R}^3}\frac{1}{2}\partial_tg^{ij}(t,x)\phi_i\phi_j\mathrm{d}x\\
&-\int_{\mathbb{R}^3} \frac{\mathrm{\partial}}
{\mathrm{\partial}x_{i}}\Big(\phi_{t}g^{ij}(t,x)\phi_j\Big)\mathrm{d}x=0.\\
\end{aligned}
\end{equation}
And since the last term is always zero, by the divergence theorem,
due to the fact that~$\phi(t,x)=0$~for~$|x|> C+t$,~\eqref{eq140}
implies
 \begin{equation}\nonumber\\
\begin{aligned}
\partial_tE_1(t)\leq CE_1(t),\\
\end{aligned}
\end{equation}
which means
 \begin{equation}\nonumber\\
\begin{aligned}
E_1(t)\leq E_1(t_0)e^{C(t-t_0)},\\
\end{aligned}
\end{equation}
so~$E_1(t)$~or~$E(t)$~is bounded, as desired.\\
\indent To prove the other half of lemma 3.2, we integrate
\eqref{eq138} over~$Q_s^t$~ and arrive at the "flux identity":
\begin{equation}\nonumber\\
\begin{aligned}
&\frac{1}{2}\int_{Q(t)}\Big(\phi_t^2(t, x)+g^{ij}(t,
x)\partial_i\phi(t, x)\partial_j\phi(t, x)+\frac{\phi^6(t,
x)}{3}\Big)\mathrm{d}x+Flux_1(\phi, M_s^t)\\
&-\frac{1}{2}\int_{Q(s)}\Big(\phi_t^2(s, x)+g^{ij}(s,
x)\partial_i\phi(s, x)\partial_j\phi(s, x)+\frac{\phi^6(s,
x)}{3}\Big)\mathrm{d}x\\
&=\frac{1}{2}\int_{Q_s^t}\partial_tg^{ij}(\tau,x)\partial_i\phi(\tau, x)
\partial_j\phi(\tau, x)\mathrm{d}x\mathrm{d}\tau,\\
\end{aligned}
\end{equation}
that is
\begin{equation}\label{eq330}
\begin{aligned}
E_1(\phi, Q(t))+Flux_1(\phi, M_s^t)-E_1(\phi, Q(s))\leq C(t_0)\int_s^t E_1(\phi, Q(\tau))\mathrm{d}\tau,\\
\end{aligned}
\end{equation}
where~$C(t_0)$~is a constant depending on~$t_0$.~And it means
\begin{equation}\label{eq331}
\begin{aligned}
&E_1(\phi, Q(t))-C(t_0)\int_{t_0}^tE_1(\phi,
Q(\tau))\mathrm{d}\tau+Flux_1(\phi,
M_s^t)\\
&\leq E_1(\phi, Q(s))-C(t_0)\int_{t_0}^sE_1(\phi,
Q(\tau))\mathrm{d}\tau,\\
\end{aligned}
\end{equation}
which implies~$E_1(\phi, Q(t))-C(t_0)\int_{t_0}^tE_1(\phi,
Q(\tau))\mathrm{d}\tau$~is a non-increasing function on~$[t_0,
0)$.~It is also bounded as we have showed above, hence~$E_1(\phi,
Q(t))-C(T)\int_{t_0}^tE_1(\phi, Q(\tau))\mathrm{d}\tau$~\\
and~$E_1(\phi, Q(s))-C(T)\int_{t_0}^sE_1(\phi,
Q(\tau))\mathrm{d}\tau$~in \eqref{eq330} must approach a
common limit. This in turn gives the important fact that\\
\begin{equation}\nonumber\\
Flux_1(\phi, M_s^t)\rightarrow 0,~~when~s, t\rightarrow 0,\\
 \end{equation}
thanks to lemma 3.1, we complete the proof of lemma 3.2.\\
\indent To prove lemma 1.2, we need to introduce the energy-momentum
tensor~$\Pi$~as a symmetric 2-tensor by
\begin{equation}\nonumber\\
\begin{aligned}
\Pi(X, Y)&=X(\phi)Y(\phi)-\frac{1}{2}<X, Y>|\nabla\phi|^2,\\
\Pi_{\alpha\beta}&=\partial_\alpha \phi\partial_\beta
\phi-\frac{1}{2}g_{\alpha\beta}|\nabla\phi|^2,\\
\end{aligned}
\end{equation}
where~$X, Y$~are vector fields and~$\phi$~ a fixed~$C^1$~function.
Then we have
\begin{equation}\label{eq144}
\begin{aligned}
\Pi(\underline{L},
\underline{L})&=\big(\underline{L}(\phi)\big)^2,~~
\Pi(L, L)=\big(L(\phi)\big)^2,\\
\Pi(\underline{L}, e_a)&=\underline{L}(\phi)e_a(\phi), ~\Pi(L,
e_a)=L(\phi)e_a(\phi),\\
\Pi(\underline{L}, L)&=L(\phi)\underline{L}(\phi)-\frac{1}{2}<L,
\underline{L}>|\nabla \phi|^2=L(\phi)\underline{L}(\phi)+|\nabla
\phi|^2=|\overline{\nabla} \phi|^2,\\
\Pi(e_a, e_b)&=e_a(\phi)e_b(\phi)-\frac{1}{2}<e_a, e_b>|\nabla
\phi|^2\\
&=e_a(\phi)e_b(\phi)-\frac{1}{2}\delta_{ab}(|\overline{\nabla}
\phi|^2-L(\phi)\underline{L}(\phi)),\\
\end{aligned}
\end{equation}
where~$\delta_{ab}$~denotes the Kronecker delta function.\\
\indent We also need a key formula showed as a lemma below.\\
 {\bf Lemma 3.3.} Let~$\phi$~be a~$C^1$~function and~$\Pi$~be
the associated energy-momentum tensor. Let~$X$~be a vector field,
and set~$P_\alpha=\Pi_{\alpha\beta}X^\beta$,~then
\begin{equation}\label{eq141}
\begin{aligned}
div P\equiv D_\alpha P^\alpha=\Box_g\phi
X(\phi)+\frac{1}{2}\Pi^{\alpha\beta} {^{(X)}\!\pi}_{\alpha\beta},\\
\end{aligned}
\end{equation}
where~$\Box_g$~is the wave operators associated to the given
metric~$g$~and has formula as follows:\\
\begin{equation}\label{eq142}
\begin{aligned}
&\Box_g\phi=|g|^{-1/2}\partial_{\alpha}(g^{\alpha\beta}|
g|^{1/2}\partial_{\beta}\phi)\\
&=-\partial_{tt}\phi+\partial_i\big(g^{ij}(t,
x)\phi_j\big)+\frac{1}{2}g^{ij}g^{lm}\partial
_mg_{ij}\partial_l\phi-\frac{1}{2}g^{ij}\partial_tg_{ij}\partial_t\phi,\\
\end{aligned}
\end{equation}
where~$|g|$~is the absolute value of the determinant of the
matrix~$(g_{\alpha\beta})$~and~$(g^{\alpha\beta})$~its inverse
matrix.\\
\indent For the proof one can read \cite{Alinhac}.\\
\indent We then construct a
multiplier:~$\frac{1}{2}(\underline{u}L+u\underline{L})+1$,~which is
close to the Morawetz multiplier~$t\partial_t+r\partial_r+1$,~and
setting~$Y=\frac{1}{2}(\underline{u}L+u\underline{L})$.~\\
 \indent
Following Christodoulou and Klainerman \cite{Chris}, the deformation
tensor of a given vector field~$X$~is the symmetric
2-tensor~$^{(X)}\!\pi$~defined by
\begin{equation}\nonumber\\
^{(X)}\!\pi(Y,Z)\equiv\pi(Y,Z)=<D_YX,Z>+<D_ZX,Y>.\\
 \end{equation}
In local coordinates
\begin{equation}\nonumber\\
\pi_{\alpha\beta}=D_{\alpha}X_{\beta}+D_{\beta}X_{\alpha},\\
 \end{equation}
as
\begin{equation}\nonumber\\
\begin{aligned}
\nabla u&=2\nabla t-\nabla \underline{u}=-2\partial_t+\underline{L}\\
&=-(m^{-1}L+m\underline{L})+\underline{L}=-m^{-1}L+(1-m)\underline{L},\\
\end{aligned}
\end{equation}
then we can compute the deformation tensor
of~$Y=\frac{1}{2}(\underline{u}L+u\underline{L})$~as follows
\begin{equation}\nonumber\\
\begin{aligned}
&^{(Y)}\!\pi_{\underline{L}\underline{L}}=0, ~~~~^{(Y)}\!\pi_{L
\underline{L}}=-2-\frac{2}{m}+2\underline{\omega}\underline{u},\\
&^{(Y)}\!\pi_{LL}=4(1-m)-4\underline{\omega}u,~~~
^{(Y)}\!\pi_{\underline{L}e_a}=(\eta_a-\underline{\eta}_a)\underline{u},\\
&^{(Y)}\!\pi_{Le_a}=\xi_a\underline{u}+2\underline{\eta}_au,~~~
^{(Y)}\!\pi_{e_ae_b}=\underline{\chi}_{ab}u+\chi_{ab}\underline{u}.\\
\end{aligned}
\end{equation}
Also
\begin{equation}\label{eq183}
\begin{aligned}
div Y&=g^{\alpha\beta}<D_\alpha Y, e_\beta>=g^{\alpha\beta}<D_\alpha
\frac{1}{2}(\underline{u}L+u\underline{L}), e_\beta>\\
&=\frac{1}{2}(\chi_{aa}\underline{u}+\chi_{bb}\underline{u}+
\underline{\chi}_{aa}u+\underline{\chi}_{bb}u)+1+m^{-1}-\underline{u}\underline{\omega}.\\
\end{aligned}
\end{equation}
Combining \eqref{eq1} and \eqref{eq142}, we get
\begin{equation}\label{eq143}
\begin{aligned}
\Box_g\phi=\phi^5+\frac{1}{2}g^{ij}g^{lm}\partial
_mg_{ij}\partial_l\phi-\frac{1}{2}g^{ij}\partial_tg_{ij}\partial_t\phi,\\
\end{aligned}
\end{equation}
together with \eqref{eq141}, substitueing~$X$~with~$Y$~we arrive at
\begin{equation}\label{eq170}
\begin{aligned}
&div P\equiv D_\alpha P^\alpha=\Box_g\phi
Y(\phi)+\frac{1}{2}\Pi^{\alpha\beta} {^{(Y)}\!\pi}_{\alpha\beta}\\
&=\big(\phi^5+\frac{1}{2}g^{ij}g^{lm}\partial
_mg_{ij}\partial_l\phi-\frac{1}{2}g^{ij}\partial_tg_{ij}\partial_t\phi\big)Y(\phi)+\frac{1}{2}\Pi^{\alpha\beta}
{^{(Y)}\!\pi}_{\alpha\beta}\\
&=Y(\frac{\phi^6}{6})+\big(\frac{1}{2}g^{ij}g^{lm}\partial
_mg_{ij}\partial_l\phi-\frac{1}{2}g^{ij}\partial_tg_{ij}\partial_t\phi\big)Y(\phi)+\frac{1}{2}\Pi^{\alpha\beta}
{^{(Y)}\!\pi}_{\alpha\beta}\\
 &=div(\frac{\phi^6Y}{6})-\frac{\phi^6}{6}divY+\big(\frac{1}{2}g^{ij}g^{lm}\partial
_mg_{ij}\partial_l\phi-\frac{1}{2}g^{ij}\partial_tg_{ij}\partial_t\phi\big)Y(\phi)+\frac{1}{2}\Pi^{\alpha\beta}
{^{(Y)}\!\pi}_{\alpha\beta}\\
\end{aligned}
\end{equation}
where~$P_\alpha=\Pi_{\alpha\beta}Y^\beta,$~and it means
\begin{equation}\label{eq171}
\begin{aligned}
-div(P-\frac{1}{6}\phi^6Y)&=\frac{1}{6}\phi^6div
Y-Y(\phi)\big(\frac{1}{2}g^{ij}g^{lm}\partial
_mg_{ij}\partial_l\phi-\frac{1}{2}g^{ij}\partial_tg_{ij}\partial_t\phi\big)\\
&-\frac{1}{2}\Pi^{\alpha\beta}
{^{(Y)}\!\pi}_{\alpha\beta}\mathop =\limits^{\vartriangle} \widetilde{R}(t, x).\\
\end{aligned}
\end{equation}
By \eqref{eq175} and \eqref{eq143}, we have
\begin{equation}\nonumber\\
\begin{aligned}
\Box_g(\frac{1}{2}\phi^2)&=div\big(\nabla(\frac{1}{2}\phi^2)\big)=div(\phi\nabla
\phi)=<D_\alpha \phi\nabla \phi, \partial^\alpha>\\
&=\partial_\alpha(\phi)<\nabla \phi,
\partial^\alpha>+\phi\Box_g\phi\\
&=g^{\alpha\beta}\partial_\alpha(\phi)<\nabla \phi,
\partial_\beta>+\phi\Box_g\phi\\
&=-(\partial_t\phi)^2+g^{ij}\phi_i\phi_j+\phi\Box_g\phi\\
&=|\overline{\nabla}\phi|^2-L(\phi)\underline{L}(\phi)+\phi^6+\phi(\frac{1}{2}g^{ij}g^{lm}\partial
_mg_{ij}\partial_l\phi-\frac{1}{2}g^{ij}\partial_tg_{ij}\partial_t\phi),\\
\end{aligned}
\end{equation}
so
\begin{equation}\label{eq178}
\begin{aligned}
-div(\phi\nabla
\phi)=-|\overline{\nabla}\phi|^2+L(\phi)\underline{L}(\phi)-\phi^6-\phi(\frac{1}{2}g^{ij}g^{lm}\partial
_mg_{ij}\partial_l\phi-\frac{1}{2}g^{ij}\partial_tg_{ij}\partial_t\phi).\\
\end{aligned}
\end{equation}
Adding \eqref{eq171} and \eqref{eq178}, we get
\begin{equation}\label{eq179}
\begin{aligned}
&-div(P-\frac{1}{6}\phi^6Y+\phi\nabla \phi)=\widetilde{R}(t, x)\\
&-|\overline{\nabla}\phi|^2+L(\phi)\underline{L}(\phi)-\phi^6-\phi(\frac{1}{2}g^{ij}g^{lm}\partial
_mg_{ij}\partial_l\phi-\frac{1}{2}g^{ij}\partial_tg_{ij}\partial_t\phi)\mathop
=\limits^{\vartriangle} R(t, x).\\
\end{aligned}
\end{equation}
Integrating the identity \eqref{eq179} over the truncated geodesic
cone~$Q_S^T,~S< T<0$,~we arrive at
\begin{equation}\nonumber\\
\begin{aligned}
&-\int_{Q(T)}<P-\frac{1}{6}\phi^6Y+\phi\nabla \phi,
-\partial_t>\mathrm{d}v-\int_{M_S^T}\frac{<P-\frac{1}{6}\phi^6Y+\phi\nabla
\phi, \nabla
\underline{u}>}{\sqrt{(\partial_t\underline{u})^2+(g^{ij}\partial_i\underline{u})^2}}\mathrm{d}\sigma\\
&+\int_{Q(S)}<P-\frac{1}{6}\phi^6Y+\phi\nabla \phi,
-\partial_t>\mathrm{d}v=\int_{Q_S^T}R(t, x)\mathrm{d}v\mathrm{d}t,\\
\end{aligned}
\end{equation}
that is
\begin{equation}\label{eq149}
\begin{aligned}
&\int_{Q(T)}\Pi(Y, \partial_t)-<\frac{1}{6}\phi^6Y-\phi\nabla \phi,
\partial_t>\mathrm{d}v-\int_{Q(S)}\Pi(Y, \partial_t)-<\frac{1}{6}\phi^6Y-\phi\nabla \phi,
\partial_t>\mathrm{d}v\\
&+\int_{M_S^T}\frac{<P-\frac{1}{6}\phi^6Y+\phi\nabla \phi,
\underline{L}>}
{\sqrt{(\partial_t\underline{u})^2+(g^{ij}\partial_i\underline{u})^2}}\mathrm{d}\sigma
=\int_{Q_S^T}R(t, x)\mathrm{d}v\mathrm{d}t.\\
\end{aligned}
\end{equation}
By \eqref{eq500}, we have
\begin{equation}\label{eq150}
\begin{aligned}
&\Pi(Y, \partial_t)-<\frac{1}{6}\phi^6Y-\phi\nabla \phi,
\partial_t>=\Pi\big(\frac{1}{2}(\underline{u}L+u\underline{L}),~
\frac{1}{2}(m^{-1}L+m\underline{L})\big)\\
&-<\frac{1}{6}\phi^6\frac{1}{2}(\underline{u}L+u\underline{L})-\phi\nabla
\phi,~
\frac{1}{2}(m^{-1}L+m\underline{L})>\\
&=\frac{mu}{4}\big(\underline{L}(\phi)\big)^2+\frac{\underline{u}}{4m}\big(L(\phi)\big)^2
+\big(\frac{u}{4m}+\frac{m\underline{u}}{4}\big)|\overline{\nabla}\phi|^2
+\big(\frac{u}{12m}+\frac{m\underline{u}}{12}\big)\phi^6\\
&+\frac{1}{2m}\phi L(\phi)+\frac{m}{2}\phi\underline{L}(\phi),\\
\end{aligned}
\end{equation}
and
\begin{equation}\label{eq151}
\begin{aligned}
&<P-\frac{1}{6}\phi^6Y+\phi\nabla \phi, ~\underline{L}>=\Pi(Y,~
\underline{L})-<\frac{1}{6}\phi^6Y-\phi\nabla \phi, ~\underline{L}>\\
&=\Pi\big(\frac{1}{2}(m^{-1}L+m\underline{L}),~
\underline{L}\big)-<\frac{1}{6}\phi^6\frac{1}{2}(\underline{u}L+u\underline{L})-\phi\nabla
\phi,~
\underline{L}>\\
&=\frac{1}{2}\underline{u}|\overline{\nabla}\phi|^2+\frac{1}{2}u\big(\underline{L}(\phi)\big)^2+
\frac{\underline{u}\phi^6}{6}+\phi\underline{L}(\phi),\\
\end{aligned}
\end{equation}
then \eqref{eq149} becomes
\begin{equation}\label{180}
\begin{aligned}
&\int_{Q(T)}\big[\frac{mu}{4}\big(\underline{L}(\phi)\big)^2+\frac{\underline{u}}{4m}\big(L(\phi)\big)^2
+\big(\frac{u}{4m}+\frac{m\underline{u}}{4}\big)|\overline{\nabla}\phi|^2
+\big(\frac{u}{12m}+\frac{m\underline{u}}{12}\big)\phi^6\\
&+\frac{1}{2m}\phi
L(\phi)+\frac{m}{2}\phi\underline{L}(\phi)\big]\mathrm{d}v+\int_{M_S^T}\frac{\frac{1}{2}u\big(\underline{L}(\phi)\big)^2
+\frac{1}{2}\underline{u}|\overline{\nabla}\phi|^2+\frac{\underline{u}}{6}\phi^6+
\phi\underline{L}(\phi)}{\sqrt{(\partial_t\underline{u})^2+(g^{ij}\partial_i\underline{u})^2}}\mathrm{d}\sigma\\
&-\int_{Q(S)}\big[\frac{mu}{4}\big(\underline{L}(\phi)\big)^2+\frac{\underline{u}}{4m}\big(L(\phi)\big)^2
+\big(\frac{u}{4m}+\frac{m\underline{u}}{4}\big)|\overline{\nabla}\phi|^2
+\big(\frac{u}{12m}+\frac{m\underline{u}}{12}\big)\phi^6\\
&+\frac{1}{2m}\phi
L(\phi)+\frac{m}{2}\phi\underline{L}(\phi)\big]\mathrm{d}v=\int_{Q_S^T}R(t,
x)\mathrm{d}v\mathrm{d}t,\\
\end{aligned}
\end{equation}
where~$Q(S)=\{x\in \mathbb{R}^3:~\underline{u}\leq 0, t=S\}$.~Noting
that~$\underline{u}=0$~on the mantle~$M_S^T$,~and when~$S, T$~is
small enough we can let~$m=1$~for the error margin is nothing
but~$\mathcal {O}(t^2)E(t)$,~then~\eqref{180} becomes a little
simpler form
\begin{equation}\label{eq181}
\begin{aligned}
&\int_{Q(T)}\big[\frac{u}{4}\big(\underline{L}(\phi)\big)^2+\frac{\underline{u}}{4}\big(L(\phi)\big)^2
+\frac{T}{2}|\overline{\nabla}\phi|^2
+\frac{T}{6}\phi^6+\frac{1}{2}\phi
L(\phi)+\frac{1}{2}\phi\underline{L}(\phi)\big]\mathrm{d}v\\
&+\int_{M_S^T}\frac{t\big(\underline{L}(\phi)\big)^2
+\phi\underline{L}(\phi)}{\sqrt{(\partial_t\underline{u})^2+(g^{ij}\partial_i\underline{u})^2}}\mathrm{d}\sigma\\
&-\int_{Q(S)}\big[\frac{u}{4}\big(\underline{L}(\phi)\big)^2+\frac{\underline{u}}{4}\big(L(\phi)\big)^2
+\frac{S}{2}|\overline{\nabla}\phi|^2
+\frac{S}{6}\phi^6+\frac{1}{2}\phi
L(\phi)+\frac{1}{2}\phi\underline{L}(\phi)\big]\mathrm{d}v\\
&=\int_{Q_S^T}R(t,
x)\mathrm{d}v\mathrm{d}t.\\
\end{aligned}
\end{equation}
Denote
\begin{equation}\nonumber\\
\begin{aligned}
&I=\int_{Q(T)}\big[\frac{u}{4}\big(\underline{L}(\phi)\big)^2+\frac{\underline{u}}{4}\big(L(\phi)\big)^2
+\frac{T}{2}|\overline{\nabla}\phi|^2
+\frac{T}{6}\phi^6+\frac{1}{2}\phi
L(\phi)+\frac{1}{2}\phi\underline{L}(\phi)\big]\mathrm{d}v,\\
&II=\int_{M_S^T}\frac{t\big(\underline{L}(\phi)\big)^2
+\phi\underline{L}(\phi)}{\sqrt{(\partial_t\underline{u})^2+(g^{ij}\partial_i\underline{u})^2}}\mathrm{d}\sigma
=\int_{M_S^0}\frac{t\big(\underline{L}(\phi)\big)^2
+\phi\underline{L}(\phi)}{\sqrt{(\partial_t\underline{u})^2+(g^{ij}\partial_i\underline{u})^2}}\mathrm{d}\sigma\\
&-\int_{M_T^0}\frac{t\big(\underline{L}(\phi)\big)^2
+\phi\underline{L}(\phi)}{\sqrt{(\partial_t\underline{u})^2+(g^{ij}\partial_i\underline{u})^2}}\mathrm{d}\sigma
=II_1-II_2,\\
&III=-\int_{Q(S)}\big[\frac{u}{4}\big(\underline{L}(\phi)\big)^2+\frac{\underline{u}}{4}\big(L(\phi)\big)^2
+\frac{S}{2}|\overline{\nabla}\phi|^2
+\frac{S}{6}\phi^6+\frac{1}{2}\phi
L(\phi)+\frac{1}{2}\phi\underline{L}(\phi)\big]\mathrm{d}v,\\
\end{aligned}
\end{equation}
then~$\eqref{eq181}$~becomes
\begin{equation}\label{eq315}
\begin{aligned}
I+II_1-II_2+III=\int_{Q_S^T}R(t,
x)\mathrm{d}v\mathrm{d}t.\\
\end{aligned}
\end{equation}
 \indent Let us estimate the right-hind side of \eqref{eq315}
first.
\begin{equation}\label{eq182}
\begin{aligned}
&\Pi^{\alpha\beta}
{^{(Y)}\!\pi}_{\alpha\beta}=g^{\alpha\alpha'}g^{\beta\beta'}\Pi_{\alpha'\beta'}{^{(Y)}\!\pi}_{\alpha\beta}\\
&=(\underline{\omega}\underline{u}-1-\frac{1}{m})|\overline{\nabla}\phi|^2+(1-m-
\underline{\omega}u)\big(\underline{L}(\phi)\big)^2\\
&-\sum_{a=1}^{2}(\eta_a-\underline{\eta}_a)\underline{u}L(\phi)e_a(\phi)+
\sum_{a,b=1}^{2}(\chi_{ab}\underline{u}+\underline{\chi}_{ab}u)e_a(\phi)e_b(\phi)\\
&-\frac{1}{2}(\chi_{aa}\underline{u}+\chi_{bb}\underline{u}+
\underline{\chi}_{aa}u+\underline{\chi}_{bb}u)|\overline{\nabla}\phi|^2\\
&+\frac{1}{2}(\chi_{aa}\underline{u}+\chi_{bb}\underline{u}+
\underline{\chi}_{aa}u+\underline{\chi}_{bb}u)L(\phi)\underline{L}(\phi)\\
&-\sum_{a=1}^{2}(\xi_a\underline{u}+
2\underline{\eta}_au)\underline{L}(\phi)e_a(\phi).\\
\end{aligned}
\end{equation}
Combining \eqref{eq183} \eqref{eq171} \eqref{eq179} with
\eqref{eq182}, and set~$m=1$~(will not influence our result) we get
\begin{equation}\label{eq188}
\begin{aligned}
\int_{Q_S^T}R(t,
x)\mathrm{d}v\mathrm{d}t&=\int_{Q_S^T}\Big[\big(\frac{\frac{1}{2}(\chi_{aa}\underline{u}+\chi_{bb}\underline{u}+
\underline{\chi}_{aa}u+\underline{\chi}_{bb}u)+2-\underline{u}\underline{\omega}}{6}-\frac{2}{3}\big)\frac{\phi^6}{6}\\
&+\big(\frac{1}{4}(\chi_{aa}\underline{u}+\chi_{bb}\underline{u}+
\underline{\chi}_{aa}u+\underline{\chi}_{bb}u)-1-\frac{1}{2}(\underline{u}\underline{\omega}-2)\big)
|\overline{\nabla}\phi|^2\\
&-\frac{1}{2}\sum_{a,b=1}^{2}(\chi_{ab}\underline{u}+\underline{\chi}_{ab}u)e_a(\phi)e_b(\phi)\\
&+\big(1-\frac{1}{4}(\chi_{aa}\underline{u}+\chi_{bb}\underline{u}+
\underline{\chi}_{aa}u+\underline{\chi}_{bb}u)\big)L(\phi)\underline{L}(\phi)-\frac{1}{2}\underline{\omega}u\big(\underline{L}(\phi)\big)^2\\
&+
\frac{1}{2}\sum_{a=1}^{2}(\eta_a-\underline{\eta}_a)\underline{u}L(\phi)e_a(\phi)+
\sum_{a=1}^{2}(\xi_a\underline{u}+
2\underline{\eta}_au)\underline{L}(\phi)e_a(\phi)\\
&-\phi(\frac{1}{2}g^{ij}g^{lm}\partial
_mg_{ij}\partial_l\phi-\frac{1}{2}g^{ij}\partial_tg_{ij}\partial_t\phi)\\
&-\frac{1}{2}\big(\underline{u}
L(\phi)+u\underline{L}(\phi)\big)(\frac{1}{2}g^{ij}g^{lm}\partial
_mg_{ij}\partial_l\phi-\frac{1}{2}g^{ij}\partial_tg_{ij}\partial_t\phi)\\
&-\frac{\phi^6}{3}\Big]\mathrm{d}v\mathrm{d}t.\\
\end{aligned}
\end{equation}
\indent Also we have
\begin{equation}\label{eq192}
\begin{aligned}
&\int_{Q_S^T}\big(-\phi(\frac{1}{2}g^{ij}g^{lm}\partial
_mg_{ij}\partial_l\phi-\frac{1}{2}g^{ij}\partial_tg_{ij}\partial_t\phi)\big)\mathrm{d}v\mathrm{d}t\\
&\leq
C(T-S)\big(\int_{Q(S)}\phi^6\mathrm{d}v\big)^{\frac{1}{6}}\big(\int_{Q(S)}\mathrm{d}v\big)^{\frac{1}{3}}
\big[\big(\int_{Q(S)}(\partial_t\phi)^2\mathrm{d}v\big)^{\frac{1}{2}}+
\big(\int_{Q(S)}(\partial_j\phi)^2\mathrm{d}v\big)^{\frac{1}{2}}\big]\\
&\leq C(T-S)|S|\big(E(\phi, Q(S))\big)^{\frac{2}{3}},\\
&\int_{Q_S^T}-\frac{1}{2}\big(\underline{u}
L(\phi)+u\underline{L}(\phi)\big)(\frac{1}{2}g^{ij}g^{lm}\partial
_mg_{ij}\partial_l\phi-\frac{1}{2}g^{ij}\partial_tg_{ij}\partial_t\phi)\mathrm{d}v\mathrm{d}t\\
&\leq C|S|(T-S)\big(E(\phi, Q(S))\big).\\
\end{aligned}
\end{equation}
Combining (2.5), (2.6), (2.8), (2.9), (2.11), (2.12), \eqref{eq196},
\eqref{eq188} and \eqref{eq192}, we get
\begin{equation}\label{eq300}
\begin{aligned}
\int_{Q_S^T}R(t, x)\mathrm{d}v\mathrm{d}t\leq C|S|(T-S)\big(E(\phi,
Q(S))\big)+C(T-S)|S|\big(E(\phi, Q(S))\big)^{\frac{2}{3}}.\\
\end{aligned}
\end{equation}
\indent On the surface~$M_S^T$~where~$\underline{u}=0$,~we have
\begin{equation}\nonumber\\
\begin{aligned}
&t\big(\underline{L}(\phi)\big)^2 +\phi\underline{L}(\phi)\\
&=t\big(m^{-1}\partial_t\phi-g^{ij}\partial_i\underline{u}\partial_j\phi\big)^2+
\phi\big(m^{-1}\partial_t\phi-g^{ij}\partial_i\underline{u}\partial_j\phi\big)\\
&=-(\underline{u}-t)\big(g^{ij}\partial_i\underline{u}\partial_j\phi-m^{-1}\partial_t\phi\big)^2
-\phi\big(g^{ij}\partial_i\underline{u}\partial_j\phi-m^{-1}\partial_t\phi\big).\\
\end{aligned}
\end{equation}
If we parameterize~$M_S^0$~by
\begin{equation}\nonumber\\
y \rightarrow \big(f(y), y\big),~~y\in Q(S),
\end{equation}
then by~$\underline{u}\big(f(y), y\big)=0$~on~$M_S^0$,~we have
\begin{equation}\nonumber\\
\begin{aligned}
&\underline{u}_tf_i+\underline{u}_i=0,\\
&f_i=-\frac{\underline{u}_i}{\underline{u}_t}=-m\underline{u}_i,\\
\end{aligned}
\end{equation}
and let~$\psi(y)=\phi\big(f(y),
y\big)$,~then~$d\sigma=\sqrt{(\partial_t\underline{u})^2+(g^{ij}\partial_i\underline{u})^2}dy$~and
\begin{equation}\nonumber\\
\begin{aligned}
\psi_j=\phi_tf_j+\phi_j,\\
\end{aligned}
\end{equation}
which implies
\begin{equation}\nonumber\\
\begin{aligned}
g^{ij}\partial_i\underline{u}\partial_j\psi&=\phi_tg^{ij}\partial_i\underline{u}
f_j+g^{ij}\partial_i\underline{u}\partial_j\phi\\
&=-m\phi_tg^{ij}\partial_i\underline{u}
\partial_j\underline{u}+g^{ij}\partial_i\underline{u}\partial_j\phi\\
&=-m\phi_t(\partial_t\phi)^2+g^{ij}\partial_i\underline{u}\partial_j\phi\\
&=-m^{-1}\phi_t+g^{ij}\partial_i\underline{u}\partial_j\phi.\\
\end{aligned}
\end{equation}
Thus, a calculation gives
\begin{equation}\label{eq306}
\begin{aligned}
II_1&=-\int_{Q(S)}\big[(\underline{u}-S)(g^{ij}\partial_i\underline{u}\partial_j\psi)^2+\psi
g^{ij}\partial_i\underline{u}\partial_j\psi \big]\mathrm{d}v\\
&=-\int_{Q(S)}\frac{\big((\underline{u}-S)g^{ij}\partial_i\underline{u}\partial_j\psi+\psi\big)^2}
{\underline{u}-S}\mathrm{d}v+\int_{Q(S)}\frac{\psi^2}{\underline{u}-S}+\psi
g^{ij}\partial_i\underline{u}\partial_j\psi \mathrm{d}v.\\
\end{aligned}
\end{equation}
Integrating by parts we see
\begin{equation}\label{eq307}
\begin{aligned}
&\int_{Q(S)}\psi g^{ij}\partial_i\underline{u}\partial_j\psi
\mathrm{d}v\\
=&\int_{Q(S)}g^{ij}\partial_i\underline{u}\partial_j(\frac{1}{2}\psi^2)\mathrm{d}v\\
=&\int_{Q(S)}\big[\partial_j(\frac{1}{2}g^{ij}\partial_i\underline{u}\psi^2)-\frac{1}{2}\psi^2\partial_j
(g^{ij}\partial_i\underline{u})\big]\mathrm{d}v\\
=&\int_{\partial
Q(S)}\frac{g^{ij}\partial_i\underline{u}\partial_j\underline{u}}{2\sqrt{\sum_{j=1}^3(
\partial_j\underline{u})^2}}\psi^2\mathrm{d}\sigma-\int_{Q(S)}\frac{1}{2}\psi^2\partial_j
(g^{ij}\partial_i\underline{u})\mathrm{d}v.\\
\end{aligned}
\end{equation}
Note that
\begin{equation}\nonumber\\
\begin{aligned}
\Box_g\underline{u}&=div(\nabla
\underline{u})=-div(\underline{L})=-\underline{\chi}_{11}-\underline{\chi}_{22}\\
&=-\partial_{tt}\underline{u}+\partial_j\big(g^{ij}(t,
x)\partial_i\underline{u}\big)+\frac{1}{2}g^{ij}g^{lm}\partial
_mg_{ij}\partial_l\underline{u}-\frac{1}{2}g^{ij}\partial_tg_{ij}\partial_t\underline{u},\\
\end{aligned}
\end{equation}
which yields\\
\begin{equation}\nonumber\\
\begin{aligned}
&\partial_j\big(g^{ij}(t,
x)\partial_i\underline{u}\big)=-\underline{\chi}_{11}-\underline{\chi}_{22}+
\partial_{tt}\underline{u}-\frac{1}{2}g^{ij}g^{lm}\partial
_mg_{ij}\partial_l\underline{u}+\frac{1}{2}g^{ij}\partial_tg_{ij}\partial_t\underline{u},\\
\end{aligned}
\end{equation}
then from (2.6) and (2.7), we have
\begin{equation}\label{eq308}
\begin{aligned}
\frac{2}{\underline{u}-t}+Ct+C\leq \partial_j(g^{ij}\partial_i\underline{u})\leq \frac{2}{\underline{u}-t}-Ct+C.\\
\end{aligned}
\end{equation}
Combining \eqref{eq306}, \eqref{eq307} and \eqref{eq308} we get
\begin{equation}\label{eq310}
\begin{aligned}
II_1&=-\int_{Q(S)}\frac{\big((\underline{u}-S)g^{ij}\partial_i\underline{u}\partial_j\psi+\psi\big)^2}
{\underline{u}-S}\mathrm{d}v+\int_{\partial
Q(S)}\frac{g^{ij}\partial_i\underline{u}\partial_j\underline{u}}{2\sqrt{\sum_{j=1}^3(
\partial_j\underline{u})^2}}\psi^2\mathrm{d}\nu\\
&-(CS+C)\int_{Q(S)}\psi^2\mathrm{d}v\\
&=-\int_{M_S^0}\frac{(\underline{u}-S)\big(-m^{-1}\phi_t+g^{ij}\partial_i\underline{u}\partial_j\phi
+\frac{\phi}{\underline{u}-S}\big)^2+(CS+C)\phi^2}
{\sqrt{(\partial_t\underline{u})^2+(g^{ij}\partial_i\underline{u})^2}}\mathrm{d}\sigma\\
&+\int_{\partial
Q(S)}\frac{g^{ij}\partial_i\underline{u}\partial_j\underline{u}}{2\sqrt{\sum_{j=1}^3(
\partial_j\underline{u})^2}}\psi^2\mathrm{d}\nu\\
&=\int_{M_S^0}\frac{S\big(\underline{L}(\phi)
+\frac{\phi}{S}\big)^2+(CS+C)\phi^2}
{\sqrt{(\partial_t\underline{u})^2+(g^{ij}\partial_i\underline{u})^2}}\mathrm{d}\sigma
+\int_{\partial
Q(S)}\frac{g^{ij}\partial_i\underline{u}\partial_j\underline{u}}{2\sqrt{\sum_{j=1}^3(
\partial_j\underline{u})^2}}\phi^2\mathrm{d}\sigma\\
&\leq
C|S|\int_{M_S^0}\big(\underline{L}(\phi))^2\mathrm{d}\sigma+C\int_{M_S^0}(\frac{1}{|S|}
+1+|S|)\phi^2\mathrm{d}\sigma\\
&+\int_{\partial
Q(S)}\frac{g^{ij}\partial_i\underline{u}\partial_j\underline{u}}{2\sqrt{\sum_{j=1}^3(
\partial_j\underline{u})^2}}\phi^2\mathrm{d}\nu\\
&\leq C|S|Flux(\phi,
M_S^0)+C(|S|+|S|^2+|S|^3)\big(\int_{M_S^0}\phi^6\mathrm{d}\sigma\big)^{\frac{1}{3}}\\
&+\int_{\partial
Q(S)}\frac{g^{ij}\partial_i\underline{u}\partial_j\underline{u}}{2\sqrt{\sum_{j=1}^3(
\partial_j\underline{u})^2}}\psi^2\mathrm{d}\nu\\
&\leq C|S|(Flux(\phi, M_S^0)+Flux(\phi,
M_S^0)^{\frac{1}{3}})+\int_{\partial
Q(S)}\frac{g^{ij}\partial_i\underline{u}\partial_j\underline{u}}{2\sqrt{\sum_{j=1}^3(
\partial_j\underline{u})^2}}\phi^2\mathrm{d}\nu.\\
\end{aligned}
\end{equation}
For~$III$,~a computation gives (also let~$m=1$~)
\begin{equation}\nonumber\\
\begin{aligned}
&\frac{u}{4}\big(\underline{L}(\phi)\big)^2+\frac{\underline{u}}{4}\big(L(\phi)\big)^2
+\frac{S}{2}|\overline{\nabla}\phi|^2
+\frac{S}{6}\phi^6+\frac{1}{2}\phi
L(\phi)+\frac{1}{2}\phi\underline{L}(\phi)\\
&=\frac{2S-\underline{u}}{4}\big(m^{-1}\partial_t\phi-g^{ij}\partial_i\underline{u}\partial_j\phi\big)^2
+\frac{\underline{u}}{4}\big(m\partial_t\phi+m^2g^{ij}\partial_i\underline{u}\partial_j\phi\big)^2\\
&+\frac{S}{2}\big[-(\partial_t\phi)^2+g^{ij}\partial_i\phi\partial_j\phi
+(m^{-1}\partial_t\phi-g^{ij}\partial_i\underline{u}\partial_j\phi\big)(m\partial_t\phi
+m^2g^{ij}\partial_i\underline{u}\partial_j\phi\big)\big]\\
&+\frac{S}{6}\phi^6+\phi\partial_t\phi\\
&=\frac{S}{2}\big(\phi_t^2+g^{ij}\partial_i\phi\partial_j\phi+\frac{\phi^6}{3}\big)+\phi_t\big(\phi+
(\underline{u}-S)g^{ij}\partial_i\underline{u}\partial_j\phi\big).\\
\end{aligned}
\end{equation}
For the second term on the right-hand side, using Cauchy-Schwartz
inequality we have
\begin{equation}\nonumber\\
\begin{aligned}
&\phi_t\big(\phi+
(\underline{u}-S)g^{ij}\partial_i\underline{u}\partial_j\phi\big)\\
&\leq |S|\big[\frac{\phi_t^2}{2}+\frac{\big(\phi+
(\underline{u}-S)g^{ij}\partial_i\underline{u}\partial_j\phi\big)^2}{2|S|^2}\big]\\
&\leq |S|\big[\frac{\phi_t^2}{2}+\frac{\big(\phi+
(\underline{u}-S)g^{ij}\partial_i\underline{u}\partial_j\phi\big)^2}{2(\underline{u}-S)^2}\big]\\
&=|S|\frac{\phi_t^2}{2}+\frac{|S|}{2}\big[\frac{\phi^2}{(\underline{u}-S)^2}
+(g^{ij}\partial_i\underline{u}\partial_j\phi)^2+\frac{2\phi
g^{ij}\partial_i\underline{u}\partial_j\phi}{\underline{u}-S}\big]\\
&\leq|S|\frac{\phi_t^2}{2}+\frac{|S|}{2}\big[\frac{\phi^2}{(\underline{u}-S)^2}
+g^{ij}\partial_i\underline{u}\partial_j\underline{u}g^{ij}\partial_i\phi\partial_j\phi+\frac{2\phi
g^{ij}\partial_i\underline{u}\partial_j\phi}{\underline{u}-S}\big]\\
&\leq|S|\frac{\phi_t^2}{2}+\frac{|S|}{2}\big[\frac{\phi^2}{(\underline{u}-S)^2}
+m^{-2}g^{ij}\partial_i\phi\partial_j\phi+\frac{2\phi
g^{ij}\partial_i\underline{u}\partial_j\phi}{\underline{u}-S}\big].\\
\end{aligned}
\end{equation}
As~$S<0$,~we get
\begin{equation}\nonumber\\
\begin{aligned}
&\frac{S}{2}\big(\phi_t^2+g^{ij}\partial_i\phi\partial_j\phi+\frac{\phi^6}{3}\big)+\phi_t\big(\phi+
(\underline{u}-S)g^{ij}\partial_i\underline{u}\partial_j\phi\big)\\
&\leq
\frac{S\phi^6}{6}-\frac{S\phi^2}{2(\underline{u}-S)^2}-\frac{S\phi
g^{ij}\partial_i\underline{u}\partial_j\phi}{\underline{u}-S},\\
\end{aligned}
\end{equation}
so
\begin{equation}\label{eq318}
\begin{aligned}
III&=-\int_{Q(S)}\big[\frac{S}{2}\big(\phi_t^2+g^{ij}\partial_i\phi\partial_j\phi+\frac{\phi^6}{3}\big)+\phi_t\big(\phi+
(\underline{u}-S)g^{ij}\partial_i\underline{u}\partial_j\phi\big)\big]\mathrm{d}v\\
&\geq
|S|\int_{Q(S)}\frac{\phi^6}{6}\mathrm{d}v+S\Big(\frac{1}{2}\int_{Q(S)}\frac{\phi^2}{(\underline{u}-S)^2}\mathrm{d}v
+\int_{Q(S)}\frac{\phi
g^{ij}\partial_i\underline{u}\partial_j\phi}{\underline{u}-S}\mathrm{d}v\Big).\\
\end{aligned}
\end{equation}
Together with \eqref{eq308}, a similar computation gives
\begin{equation}\label{eq311}
\begin{aligned}
&\int_{Q(S)}\frac{\phi
g^{ij}\partial_i\underline{u}\partial_j\phi}{\underline{u}-S}\mathrm{d}v\\
&=\int_{Q(S)}\frac{
g^{ij}\partial_i\underline{u}\partial_j(\frac{\phi^2}{2})}{\underline{u}-S}\mathrm{d}v\\
&=\int_{Q(S)}\partial_j\Big(\frac{
g^{ij}\partial_i\underline{u}(\frac{\phi^2}{2})}{\underline{u}-S}\Big)\mathrm{d}v-
\int_{Q(S)}\frac{\phi^2}{2}\partial_j\Big(\frac{
g^{ij}\partial_i\underline{u}}{\underline{u}-S}\Big)\mathrm{d}v\\
&=\int_{\partial
Q(S)}\frac{g^{ij}\partial_i\underline{u}\partial_j\underline{u}\phi^2}{2(\underline{u}-S)\sqrt{\sum_{j=1}^3(
\partial_j\underline{u}})^2}\mathrm{d}\nu-\int_{Q(S)}\frac{\phi^2}{2(\underline{u}-S)}
\partial_j(g^{ij}\partial_i\underline{u})\mathrm{d}v\\
&+\int_{Q(S)}\frac{\phi^2g^{ij}\partial_i\underline{u}
\partial_j\underline{u}}{2(\underline{u}-S)^2}\mathrm{d}v\\
&=\int_{\partial
Q(S)}\frac{g^{ij}\partial_i\underline{u}\partial_j\underline{u}\phi^2}{-2S\sqrt{\sum_{j=1}^3(
\partial_j\underline{u}})^2}\mathrm{d}\nu-\int_{Q(S)}\frac{\phi^2}{2(\underline{u}-S)}
(\frac{2}{\underline{u}-S}+CS+C)\mathrm{d}v\\
&+\int_{Q(S)}\frac{m^{-2}\phi^2}{2(\underline{u}-S)^2}\mathrm{d}v.\\
\end{aligned}
\end{equation}
Combining \eqref{eq318} and \eqref{eq311}, we get
\begin{equation}\label{eq312}
\begin{aligned}
III\geq |S|\int_{Q(S)}\frac{\phi^6}{6}\mathrm{d}v-\int_{\partial
Q(S)}\frac{g^{ij}\partial_i\underline{u}\partial_j\underline{u}\phi^2}{2\sqrt{\sum_{j=1}^3(
\partial_j\underline{u}})^2}\mathrm{d}\sigma-(CS+CS^2)\int_{Q(S)}\frac{\phi^2}{2(\underline{u}-S)}
\mathrm{d}v.\\
\end{aligned}
\end{equation}
Using H\"{o}lder's inequality it is easy to see that
\begin{equation}\label{eq313}
\begin{aligned}
I&=\int_{Q(T)}\big[\frac{u}{4}\big(\underline{L}(\phi)\big)^2+\frac{\underline{u}}{4}\big(L(\phi)\big)^2
+\frac{T}{2}|\overline{\nabla}\phi|^2
+\frac{T}{6}\phi^6+\frac{1}{2}\phi
L(\phi)+\frac{1}{2}\phi\underline{L}(\phi)\big]\mathrm{d}v\\
&\leq C|T|E(\phi,
Q(T))+C|T|\Big[\Big(\int_{Q(T)}\phi^6\mathrm{d}v\Big)^{\frac{1}{6}}
\Big(\big(\int_{Q(T)}(L(\phi))^2\mathrm{d}v\big)^{\frac{1}{2}}+
\big(\underline{L}(\phi)^2\big)^{\frac{1}{2}}\mathrm{d}v\Big)\Big]\\
&\leq C|T|E(\phi, Q(T))+C|T|E(\phi, Q(T))^{\frac{2}{3}},\\
II_2&=\int_{M_T^0}\frac{t\big(\underline{L}(\phi)\big)^2
+\phi\underline{L}(\phi)}{\sqrt{(\partial_t\underline{u})^2+(g^{ij}\partial_i\underline{u})^2}}\mathrm{d}\sigma
\leq C|T|Flux(\phi, M_T^0)+C|T|Flux(\phi, M_T^0)^{\frac{2}{3}}.\\
\end{aligned}
\end{equation}
Now, we combine \eqref{eq315}, \eqref{eq300}, \eqref{eq310},
\eqref{eq312} and \eqref{eq313} to obtain
\begin{equation}\nonumber\\
\begin{aligned}
|S|\int_{Q(S)}\frac{\phi^6}{6}\mathrm{d}v&\leq III+\int_{\partial
Q(S)}\frac{g^{ij}\partial_i\underline{u}\partial_j\underline{u}\phi^2}{2\sum_{j=1}^3(
\partial_j\underline{u})^2}\mathrm{d}\sigma+(CS+CS^2)\int_{Q(S)}\frac{\phi^2}{2(\underline{u}-S)}\mathrm{d}v\\
&=-I-II_1+II_2+\int_{Q_S^T}R(t,
x)\mathrm{d}v\mathrm{d}t+\int_{\partial
Q(S)}\frac{g^{ij}\partial_i\underline{u}\partial_j\underline{u}\phi^2}{2\sum_{j=1}^3(
\partial_j\underline{u})^2}\mathrm{d}\sigma\\
&+(CS+CS^2)\int_{Q(S)}\frac{\phi^2}{2(\underline{u}-S)}\mathrm{d}v\\
&\leq C|T|\big(E(\phi, Q(T)+E(\phi,
Q(T))^{\frac{2}{3}}\big)\\
&+C|S|\big(Flux(\phi, M_S^0)+Flux(\phi, M_S^0)^{\frac{1}{3}}\big)\\
&+C|T|\big(Flux(\phi, M_T^0)+Flux(\phi, M_T^0)^{\frac{2}{3}}\big)\\
&+C|S|(T-S)\big(E(\phi, Q(S))+\big(E(\phi,
Q(S))\big)^{\frac{2}{3}}\big)\\
&+(CS^2+CS^3)\big(E(\phi, Q(S))\big)^{\frac{1}{3}},\\
\end{aligned}
\end{equation}
and then the result of lemma 1.2 follows as we can
choose~$T=-S^2$.~\\
\section*{\textbf{Acknowledgement}} \indent We are very grateful
to Professor Alinhac for giving a series of lectures in Fudan
University, introducing the idea of null frame by Christodoulou and
Klainerman to us and for many helpful discussions. Also we thank
Professor Yuxin Dong and Professor Yuanlong Xin for helping us to
understand some knowledge of
Riemannian Geometry.\\
\indent The authors are supported by the National Natural Science
Foundation of China under grant 10728101, the 973 Project of the
Ministry of Science and Technology of China, the doctoral program
foundation of the Ministry Education of China, the "111" project and
SGST 09DZ2272900, the outstanding doctoral science foundation
program of Fudan University.

\end{document}